\documentclass{amsart}
\usepackage{amssymb}
\usepackage{amsmath}
\usepackage{amsfonts}

\newtheorem{theorem}{Theorem}
\theoremstyle{plain}

\newtheorem{definition}[theorem]{Definition}

\newtheorem{lemma}[theorem]{Lemma}
\newtheorem{notation}[theorem]{Notation}

\newtheorem{proposition}[theorem]{Proposition}
\newtheorem{remark}[theorem]{Remark}

\numberwithin{equation}{section}
\numberwithin{theorem}{subsection}
\input{tcilatex}

\begin{document}
\title[Quantization of symplectic vector spaces over finite fields ]{%
Quantization of symplectic vector spaces over finite fields}
\author{Shamgar Gurevich}
\address{Department of Mathematics, University of California, Berkeley, CA
94720, USA. }
\email{shamgar@math.berkeley.edu}
\author{Ronny Hadani}
\address{Department of Mathematics, University of Chicago, IL 60637, USA.}
\email{hadani@math.uchicago.edu}

\begin{abstract}
In this paper, we construct a quantization functor, associating a complex
vector space $\mathcal{H}(V)$ to a finite dimensional symplectic vector
space $V$ over a finite field of odd characteristic. As a result, we obtain
a canonical model for the Weil representation of the symplectic group $%
Sp\left( V\right) $. The main new technical result is a proof of a stronger
form of the Stone--von Neumann property for the Heisenberg group $H(V)$. Our
result answers, for the case of the Heisenberg group, a question of Kazhdan
about the possible existence of a canonical vector space attached to a
coadjoint orbit of a general unipotent group over finite field.
\end{abstract}

\maketitle

\section{Introduction}

Quantization is a fundamental procedure in mathematics and in physics.
Although it is widely used in both contexts, its precise nature remains to
some extent unclear. From the physical side, quantization is the procedure
by which one associates to a classical mechanical system its quantum
counterpart. From the mathematical side, it seems that quantization is a way
to construct interesting vector spaces out of symplectic manifolds,
suggesting a method for constructing representations of the corresponding
groups of symplectomorphisms \cite{Ki, Ko}.

\subsection{Main results}

\subsubsection{Quantization functor}

In this paper, we construct a quantization functor 
\begin{equation*}
\mathcal{H}:\mathsf{Symp}\rightarrow \mathsf{Vect},
\end{equation*}%
where $\mathsf{Symp}$ denotes the (groupoid) category whose objects are
finite dimensional symplectic vector spaces over the finite field $\mathbb{F}%
_{q}$, and morphisms are linear isomorphisms of symplectic vector spaces and 
$\mathsf{Vect}$ denotes the category of finite dimensional complex vector
spaces.

\subsubsection{Canonical model for the Weil representation}

As a consequence, for a fixed symplectic vector space $V\in $ $\mathsf{Symp}$%
, we obtain, by functoriality, a homomorphism 
\begin{equation*}
\mathcal{H}:Sp\left( V\right) \rightarrow GL\left( \mathcal{H}\left(
V\right) \right) ,
\end{equation*}%
which is isomorphic to the Weil representation \cite{W} of the group $%
Sp\left( V\right) $ and we refer to it as the \textit{canonical model} of
the Weil representation.

\begin{remark}
The symplectic group acts on the canonical model $\mathcal{H}\left( V\right) 
$ in a natural way. This should be contrasted with the classical approach 
\cite{Ge, Ho, W} in which one constructs a projective representation of the
symplectic group (invoking the Stone--von Neumann theorem) and then deriving
the linear action, implicitly, by invoking the fact that every projective
representation of $Sp\left( V\right) $ can be linearized.
\end{remark}

\subsubsection{Properties of the quantization functor}

We show that the functor $\mathcal{H}$ is a monoidal functor, compatible
with duality and with the operation of linear symplectic reduction. The last
compatibility means that given $V\in \mathsf{Symp}$ and a pair $I^{\circ
}=\left( I,o_{I}\right) $, where $I\subset V$ is an isotropic subspace and $%
o_{I}\in \wedge ^{top}I$ is a non-zero vector\footnote{%
We will call the pair $I^{\circ }=\left( I,o_{I}\right) $ an \underline{%
oriented isotropic subspace}.}, there exists a natural isomorphism 
\begin{equation}
\mathcal{H}^{I}\left( V\right) \simeq \mathcal{H}\left( I^{\bot }/I\right) ,
\label{reduc_eq}
\end{equation}%
where $\mathcal{H}^{I}\left( V\right) $ stands for the subspace of $I$%
-invariant vectors in $\mathcal{H}\left( V\right) $ \ (an operation which
will be made precise in the sequel) and $I^{\perp }/I\in \mathsf{Symp}$ is
the symplectic reduction of $V$ with respect to $I$ \cite{BW}. A particular
situation is when $I=L$ is a Lagrangian subspace. In this situation, $%
L^{\perp }/L=0$ and (\ref{reduc_eq}) yields an isomorphism $\mathcal{H}%
^{L}\left( V\right) \simeq \mathcal{H}\left( 0\right) =%
%TCIMACRO{\U{2102} }%
%BeginExpansion
\mathbb{C}
%EndExpansion
$, which associates to $1\in 
%TCIMACRO{\U{2102} }%
%BeginExpansion
\mathbb{C}
%EndExpansion
$ a well defined vector in $\mathcal{H}\left( V\right) $. This establishes a
mechanism which associates to every pair $L^{\circ }=\left( L,o_{L}\right) $
a well defined vector $v_{L^{\circ }}\in \mathcal{H}_{V}$. To the best of
our knowledge (cf. \cite{GS, We1, We2}), this kind of structure, which
exists in the setting of the Weil representation of the group $Sp\left(
V\right) ,$ was not observed before.

\subsubsection{The strong Stone--von Neumann theorem}

The main technical result of this paper is a proof (\cite{GH, H}
unpublished) of a strong form of the Stone--von Neumann theorem for the
Heisenberg group over $\mathbb{F}_{q}$ (see also the survey \cite{GH2}). We
supply two proofs; in the first proof we use only basic considerations from
linear algebra and in the second proof we construct an algebra geometric
object (an $\ell $-adic perverse Weil sheaf) which is interesting in its own
sake and in particular implies the strong Stone--von Neumann theorem.

We devote the rest of the introduction to an intuitive explanation of the
main ideas and results of this paper.

\subsection{The Heisenberg representation}

\subsubsection{The Heisenberg group}

Let $\left( V,\omega \right) $ be a $2n$-dimensional symplectic vector space
over the finite field $\mathbb{F}_{q}$. The vector space $V$ considered as
an abelian group admits a non-trivial central extension called the \textit{%
Heisenberg group} 
\begin{equation*}
0\rightarrow \mathbb{F}_{q}\rightarrow H\left( V\right) \rightarrow
V\rightarrow 0\text{. }
\end{equation*}%
\ 

Choosing a section $V\rightarrow H\left( V\right) $, the Heisenberg group
can be presented as $H\left( V\right) =V\times \mathbb{F}_{q}$ with the
center given by $Z=Z(H)=\left \{ \left( 0,z\right) :z\in \mathbb{F}%
_{q}\right \} $. The symplectic group $Sp\left( V\right) $ acts on $H\left(
V\right) $ by group automorphisms via its tautological action on the $V$%
-coordinate.

\subsubsection{The Heisenberg representation}

The representation theory of the Heisenberg group is relatively simple.
Given a non-trivial central character $\psi :Z\rightarrow 
%TCIMACRO{\U{2102} }%
%BeginExpansion
\mathbb{C}
%EndExpansion
^{\times }$, there exist a unique (up to isomorphism) irreducible
representation 
\begin{equation*}
\pi :H\left( V\right) \rightarrow GL(\mathcal{H})\text{,}
\end{equation*}%
such that the center acts by $\psi $, i.e., $\pi _{|Z}=\psi \cdot Id_{%
\mathcal{H}}$. This is the content of the celebrated Stone--von Neumann
Theorem. The representation $\pi $ is referred to as the \textit{Heisenberg
representation. }

\subsubsection{Models of the Heisenberg representation}

The Heisenberg representation admits a family of models which are associated
with\  \underline{oriented Lagrangian}\textit{\ }subspaces in $V$ \cite{B, LV}%
. An oriented Lagrangian subspace is a pair $L^{\circ }=\left(
L,o_{L}\right) $, where $L\subset V$ is a Lagrangian subspace and $o_{L}\in
\wedge ^{top}L$ is a non-zero vector. Every oriented Lagrangian $L^{\circ }$
is associated with a model 
\begin{equation*}
\pi _{L^{\circ }}:H\left( V\right) \rightarrow GL\left( \mathcal{H}%
_{L^{\circ }}\right) \text{,}
\end{equation*}%
where the vector space $\mathcal{H}_{L^{\circ }}$ is the space of functions $%
%TCIMACRO{\U{2102} }%
%BeginExpansion
\mathbb{C}
%EndExpansion
\left( L\backslash H\left( V\right) ,\psi \right) $, consisting of functions 
$f:H\left( V\right) \rightarrow 
%TCIMACRO{\U{2102} }%
%BeginExpansion
\mathbb{C}
%EndExpansion
$ such that $f\left( z\cdot l\cdot h\right) =\psi \left( z\right) f\left(
h\right) $, for every $z\in Z$, $l\in L$ and the action $\pi _{L^{\circ }}$
is given by right translation.

\begin{remark}
The definition of the model $\left( \pi _{L^{\circ }},H\left( V\right) ,%
\mathcal{H}_{L^{\circ }}\right) $ does not depend on the orientation $o_{L}$%
. The role of the orientation will appear later when we explain the relation
between different models---the strong Stone--von Neumann property.
\end{remark}

\subsection{The strong Stone--von Neumann property}

The collection of models $\left \{ \mathcal{H}_{L^{\circ }}\right \} $ can
be thought of as a vector bundle $\mathfrak{H}$ on the set $OLag\left(
V\right) $ of oriented Lagrangians, with fibers $\mathfrak{H}_{|L^{\circ }}=%
\mathcal{H}_{L^{\circ }}$. The vector bundle $\mathfrak{H}$ is equipped with
the following tautological structures:

\begin{itemize}
\item An action of the group $H\left( V\right) $ on each fiber.

\item An $Sp\left( V\right) $-equivariant structure.
\end{itemize}

The content of the strong Stone--von Neumann theorem is that $\mathfrak{H}$
admits a canonical trivialization:

\begin{theorem}[Strong Stone--von Neumann theorem]
There exists a canonical system of intertwining morphisms $T_{M^{\circ
},L^{\circ }}\in \mathrm{Hom}_{H\left( V\right) }\left( \mathcal{H}%
_{L^{\circ }},\mathcal{H}_{M^{\circ }}\right) $, for every $\left( M^{\circ
},L^{\circ }\right) \in OLag\left( V\right) ^{2}$, satisfying the following
multiplicativity property 
\begin{equation*}
T_{N^{\circ },M^{\circ }}\circ T_{M^{\circ },L^{\circ }}=T_{N^{\circ
},L^{\circ }},
\end{equation*}%
for every $\left( N^{\circ },M^{\circ },L^{\circ }\right) \in OLag\left(
V\right) ^{3}$.
\end{theorem}

\begin{remark}
The definition of the trivialization $\{T_{M^{\circ },L^{\circ }}\}$ depends
on the orientation structure.
\end{remark}

\subsubsection{Canonical vector space}

The vector space $\mathcal{H}\left( V\right) $ is the space of "horizontal
sections" of $\mathfrak{H}$%
\begin{equation*}
\mathcal{H}\left( V\right) =\Gamma _{hor}\left( OLag\left( V\right) ,%
\mathfrak{H}\right) ,
\end{equation*}%
where a vector in $\mathcal{H}\left( V\right) $ is a compatible system $%
\left( v_{L^{\circ }}\in \mathcal{H}_{L^{\circ }}:L^{\circ }\in OLag\left(
V\right) \right) $ such that $T_{M^{\circ },L^{\circ }}\left( v_{L^{\circ
}}\right) =v_{M^{\circ }}$, for every $\left( M^{\circ },L^{\circ }\right)
\in OLag\left( V\right) ^{2}$. The symplectic group $Sp\left( V\right) $
acts on the vector space $\mathcal{H}\left( V\right) $ in an obvious manner.
We denote this action by 
\begin{equation*}
\rho _{V}:Sp\left( V\right) \rightarrow GL\left( \mathcal{H}\left( V\right)
\right) ,
\end{equation*}
and refer to it as the \textit{canonical model }of the Weil representation.

We proceed to explain an algebra geometric construction of the
trivialization $\{T_{M^{\circ },L^{\circ }}\}$.

\subsubsection{Canonical trivialization}

The construction will be close in spirit to the procedure of
\textquotedblleft analytic continuation\textquotedblright . Let $%
U_{2}\subset OLag\left( V\right) ^{2}$ denote the subset consisting of pairs 
$(M^{\circ },L^{\circ })\in OLag\left( V\right) ^{2}$ such that $M+L=V$. On
the set $U_{2}$, the canonical intertwining morphisms are given by a uniform
explicit formula---\textbf{ansatz}. The problem is how to extend this
formula to the set of all pairs. \ An appealing way to do this is through
the use of algebraic geometry.

\subsubsection{Kernel functions}

Every intertwining morphism $T_{M^{\circ },L^{\circ }}$ can be uniquely
presented by a function $K_{M^{\circ },L^{\circ }}\in 
%TCIMACRO{\U{2102} }%
%BeginExpansion
\mathbb{C}
%EndExpansion
\left( M\backslash H\left( V\right) /L,\psi \right) $, which we refer to as
the \textit{canonical kernel function}. The collection of kernel functions $%
\left \{ K_{M^{\circ },L^{\circ }}:\left( M^{\circ },L^{\circ }\right) \in
U_{2}\right \} $ can be thought of as a single function $K_{O}$ on $%
O=U_{2}\times H\left( V\right) $ which is given by $K_{O}\left( M^{\circ
},L^{\circ },-\right) =K_{M^{\circ },L^{\circ }}\left( -\right) $, for every 
$\left( M^{\circ },L^{\circ }\right) \in U_{2}$. The problem now translates
to the problem of extending the function $K_{O}$ to the set $\ X=OLag\left(
V\right) ^{2}\times H\left( V\right) $. To this end, we invoke the procedure
of geometrization.

\subsection{Geometrization}

A general ideology due to Grothendieck is that any meaningful set-theoretic
object is governed by a more fundamental algebra-geometric one. \textit{\ }%
The procedure by which one translate from the set theoretic setting to
algebraic geometry is called \textit{geometrization}, which is a formal
procedure by which sets are replaced by algebraic varieties and functions
are replaced by certain sheaf-theoretic objects.

The precise setting consists of:

\begin{itemize}
\item A set $X=\mathbf{X}\left( \mathbb{F}_{q}\right) $ of rational points
of an algebraic variety $\mathbf{X}$, defined over $\mathbb{F}_{q}$.

\item A complex valued function $f\in 
%TCIMACRO{\U{2102} }%
%BeginExpansion
\mathbb{C}
%EndExpansion
\left( X\right) $ governed by an $\ell $-adic Weil sheaf $\mathcal{F}$.
\end{itemize}

The variety $\mathbf{X}$ is a space equipped with an automorphism $Fr:%
\mathbf{X}\rightarrow \mathbf{X}$ \ (called Frobenius), such that the set $X$
is naturally identified with the set of fixed points $X=\mathbf{X}^{Fr}$. \ 

The sheaf $\mathcal{F}$ can be thought of as a vector bundle on the variety $%
\mathbf{X}$, equipped with an endomorphism $\theta :\mathcal{F\rightarrow F}$
which lifts $Fr$.

The relation between the function $f$ and the sheaf $\mathcal{F}$ is called
Grothendieck's \textit{sheaf-to-function correspondence: }Given a point $%
x\in X$, the endomorphism $\theta $ restricts to an endomorphism $\theta
_{x}:\mathcal{F}_{|x}\rightarrow \mathcal{F}_{|x}$ of the fiber $\mathcal{F}%
_{|x}$. The value of $f$ on the point $x$ is given by 
\begin{equation*}
f(x)=f^{\mathcal{F}}\left( x\right) =Tr(\theta _{x}:\mathcal{F}%
_{|x}\rightarrow \mathcal{F}_{|x}).
\end{equation*}

\subsection{Solution to the extension problem}

The extension problem of the function $K_{O}$ fits nicely to the
geometrization setting:

\begin{itemize}
\item The sets $O,X$ are sets of rational points of corresponding algebraic
varieties $O=\mathbf{O}\left( \mathbb{F}_{q}\right) $\textbf{\ }and $X=%
\mathbf{X}\left( \mathbb{F}_{q}\right) $. In addition, the imbedding $%
j:O\hookrightarrow X$ is induced from an open imbedding $j:\mathbf{%
O\hookrightarrow X}.$

\item The function $K_{O}$ comes from a Weil sheaf $\mathcal{K}_{\mathbf{O}}$
on the variety $\mathbf{O}$%
\begin{equation*}
K_{O}=f^{\mathcal{K}_{\mathbf{O}}}\text{.}
\end{equation*}
\end{itemize}

The extension problem is solved as follows: First extend the sheaf $\mathcal{%
K}_{\mathbf{O}}$ to a sheaf $\mathcal{K}$ on the variety $\mathbf{X}$ and
then take the corresponding function $K=f^{\mathcal{K}}$, which establishes
the desired extension. The reasoning behind this strategy is that in the
realm of sheaves there exist several functorial operations of extension,
probably the most interesting one is called \textit{perverse extension} 
\textit{\cite{BBD}}. The sheaf $\mathcal{K}$ is defined as the perverse
extension of $\mathcal{K}_{\mathbf{O}}$.

\subsection{Structure of the paper}

Apart from the introduction, the paper consists of two sections and two
appendices.

In Section \ref{pre_sec}, the classical construction of the Weil
representation is described. We begin with the definition of the Heisenberg
group and the Heisenberg representation, then we briefly describe the
classical construction of the Weil representation. In Section \ref%
{canonical_sec} we develop the framework of canonical vector spaces.
Specifically, we introduce the canonical system of intertwining morphisms
between different models of the Heisenberg representation and formulate the
strong Stone--von Neumann property of the Heisenberg representation (Theorem %
\ref{SS-vN_thm}). Using Theorem \ref{SS-vN_thm}, we construct a quantization
functor $\mathcal{H}$. We finish this section by showing that $\mathcal{H}$
is a monoidal functor and it is compatible with duality and the operation of
linear symplectic reduction. In section \ref{geom_sec}, we construct a sheaf
theoretic counterpart of the canonical system of intertwining morphisms.
Finally, in Appendix \ref{proofs_sec}, we provide proofs for all statements
which appear in the body of the paper.

\subsection{Acknowledgements}

We would like to thank our scientific advisor J. Bernstein for his interest
and guidance, and for his idea about the notion of oriented Lagrangian
subspace. It is a pleasure to thank D. Kazhdan for sharing with us his
thoughts \cite{K} about the possible existence of canonical Hilbert spaces.
We thank A. Weinstein for several very influential discussions and for the
opportunity to present this work in the symplectic geometry seminar,
Berkeley, February 2007. We would like to acknowledge M. Vergne for the
encouragement to write this paper. Finally, we would like to thank O. Ceyhan
and the organizers of the conference AGAQ, held in Istanbul during June
2006, for the invitation to present this work.

\section{The Weil representation\label{pre_sec}}

\subsection{The Heisenberg group}

Let $(V,\omega )$ be a $2n$--dimensional symplectic vector space over the
finite field $\mathbb{F}_{q}$, where $q$ is odd. Considering $V$ as an
abelian group, it admits a non--trivial central extension $H\left( V\right) $
called the \textit{Heisenberg }group, namely 
\begin{equation*}
0\rightarrow \mathbb{F}_{q}\rightarrow H\left( V\right) \rightarrow
V\rightarrow 0.
\end{equation*}

Concretely, the group $H\left( V\right) $ can be presented as the set $%
H\left( V\right) =V\times \mathbb{F}_{q}$ with the multiplication given by%
\begin{equation*}
(v,z)\cdot (v^{\prime },z^{\prime })=(v+v^{\prime },z+z^{\prime }+\tfrac{1}{2%
}\omega (v,v^{\prime })).
\end{equation*}

The center of $H\left( V\right) $ is $\ Z=Z_{H\left( V\right) }=\{(0,z):$ $%
z\in \mathbb{F}_{q}\}.$ The symplectic group $Sp\left( V\right) =Sp(V,\omega
)$ acts by automorphism on $H\left( V\right) $ through its tautological
action on the $V$--coordinate.

\subsection{The Heisenberg representation \label{HR}}

One of the most important attributes of the group $H\left( V\right) $ is
that it admits, principally, a unique irreducible representation---this is
the Stone--von Neumann property (S-vN for short). The precise statement goes
as follows. Let $\psi :Z\rightarrow 
%TCIMACRO{\U{2102} }%
%BeginExpansion
\mathbb{C}
%EndExpansion
^{\times }$ be a non-trivial character of the center. For example we can
take $\psi \left( z\right) =e^{\frac{2\pi i}{p}tr\left( z\right) }$.

\begin{theorem}[Stone--von Neuman property]
\label{S-vN_thm}There exists a unique (up to isomorphism) irreducible
unitary representation $(\pi ,H,\mathcal{H)}$ with the center acting by $%
\psi ,$ i.e., $\pi _{|Z}=\psi \cdot Id_{\mathcal{H}}$.
\end{theorem}

The representation $\pi $ which appears in the above theorem will be called
the \textit{Heisenberg representation}.

\subsection{The Weil representation}

A direct consequence of Theorem \ref{S-vN_thm} is the existence of a
projective representation $\widetilde{\rho }:Sp\left( V\right) \rightarrow
PGL(\mathcal{H)}$. The construction of $\widetilde{\rho }$ out of the
Heisenberg representation $\pi $ is due to Weil \cite{W} and it goes as
follows. Considering the Heisenberg representation $\pi $ and an element $%
g\in Sp\left( V\right) $, one can define a new representation $\pi ^{g}$
acting on the same Hilbert space via $\pi ^{g}\left( h\right) =\pi \left(
g\left( h\right) \right) $. Clearly both $\pi $ and $\pi ^{g}$ have the same
central character $\psi $ hence by Theorem \ref{S-vN_thm} they are
isomorphic. Since the space $\mathrm{Hom}_{H\left( V\right) }(\pi ,\pi ^{g})$
is one-dimensional, choosing for every $g\in Sp\left( V\right) $ a non-zero
representative $\widetilde{\rho }(g)\in \mathrm{Hom}_{H\left( V\right) }(\pi
,\pi ^{g})$ gives the required projective representation. In more concrete
terms, the projective representation $\widetilde{\rho }$ is characterized by
the formula 
\begin{equation}
\widetilde{\rho }\left( g\right) \pi \left( h\right) \widetilde{\rho }\left(
g^{-1}\right) =\pi \left( g\left( h\right) \right) ,  \label{Egorov}
\end{equation}%
for every $g\in Sp\left( V\right) $ and $h\in H\left( V\right) $.

It is known \cite{Ge, Ho} that

\begin{theorem}[The Weil representation]
\label{Weilrep_thm}There exists a linear representation $\rho :Sp\left(
V\right) \rightarrow GL\left( \mathcal{H}\right) $ lying over $\widetilde{%
\rho }$.
\end{theorem}

In the next section we will show that the linear representation $\rho $
appears as a consequence of the existence of a canonical vector space $%
\mathcal{H(}V)$ associated with the symplectic vector space $V.$

\section{Canonical vector space\label{canonical_sec}}

\subsection{Models of the Heisenberg representation}

Although the representation $\pi $ is unique, it admits a multitude of
different models (realizations), in fact this is one of its most interesting
and powerful attributes. These models appear in families, in this work we
will be interested in a particular family of such models which are
associated with Lagrangian subspaces in $V$.

Let $Lag\left( V\right) $ denote the set of Lagrangian subspaces in $V$. Let 
$%
%TCIMACRO{\U{2102} }%
%BeginExpansion
\mathbb{C}
%EndExpansion
\left( H\left( V\right) ,\psi \right) $ denote the subspace of functions $f$ 
$\in 
%TCIMACRO{\U{2102} }%
%BeginExpansion
\mathbb{C}
%EndExpansion
\left( H\left( V\right) \right) $ satisfying the equivariance property $%
f\left( z\cdot h\right) =\psi \left( z\right) f\left( h\right) $ for every $%
z\in Z$.

We associate with each Lagrangian subspace $L\in Lag\left( V\right) $ a
model $(\pi _{L},H\left( V\right) ,\mathcal{H}_{L})$ of the Heisenberg
representation as follows: The vector space $\mathcal{H}_{L}$ consists of
functions $f\in 
%TCIMACRO{\U{2102} }%
%BeginExpansion
\mathbb{C}
%EndExpansion
\left( H\left( V\right) ,\psi \right) $ satisfying $f\left( l\cdot h\right)
=f\left( h\right) $ for every $l\in L$ and the Heisenberg action is given by
right translation $\pi _{L}\left( h\right) [f]\left( h^{\prime }\right)
=f\left( h^{\prime }\cdot h\right) $ for every $f\in \mathcal{H}_{L}$.

\begin{definition}
An \underline{oriented Lagrangian} $L^{\circ }$ is a pair $L^{\circ }=\left(
L,o_{L}\right) $, where $L\in Lag\left( V\right) $ and $o_{L}\in $ $%
\bigwedge \nolimits^{top}L$ is a non-zero vector.
\end{definition}

Let us denote by $OLag\left( V\right) $ the set of oriented Lagrangian
subspaces in $V$. Similarly, we associate with each oriented Lagrangian $%
L^{\circ }\in OLag\left( V\right) $ a model $(\pi _{L^{\circ }},H,\mathcal{H}%
_{L^{\circ }})$ of the Heisenberg representation, taking $\mathcal{H}%
_{L^{\circ }}=\mathcal{H}_{L}$ and $\pi _{L^{\circ }}=\pi _{L}$. The
collection of models $\left \{ \mathcal{H}_{L^{\circ }}\right \} $ forms a
vector bundle $\mathfrak{H}\left( V\right) \rightarrow OLag\left( V\right) $
with fibers $\mathfrak{H}_{L^{\circ }}=\mathcal{H}_{L^{\circ }}$. The vector
bundle $\mathfrak{H=H}\left( V\right) $ is equipped with an additional
structure of an action $\pi _{L^{\circ }}$ of $H\left( V\right) $ on each
fiber. This suggests the following terminology:

\begin{definition}
Let $k\in 
%TCIMACRO{\U{2115} }%
%BeginExpansion
\mathbb{N}
%EndExpansion
$. An $H\left( V\right) ^{k}$-vector bundle on $OLag\left( V\right) $ is a
vector bundle $\mathfrak{E}\rightarrow OLag\left( V\right) $ equipped with a
fiberwise action $\pi _{L^{\circ }}:H\left( V\right) ^{k}\rightarrow GL(%
\mathfrak{E}_{L^{\circ }})$, for every $L^{\circ }\in OLag\left( V\right) $.
\end{definition}

In addition, our $\mathfrak{H}$ is equipped with a natural $Sp\left(
V\right) $-equivariant structure, defined as follows: For every $g\in
Sp\left( V\right) $, let $g^{\ast }\mathfrak{H}$ be the $H\left( V\right) $%
-vector bundle with fibers $g^{\ast }\mathfrak{H}_{L^{\circ }}=\mathcal{H}%
_{gL^{\circ }}$ and the $g$-twisted Heisenberg action $\pi _{L^{\circ
}}^{g}:H\left( V\right) \rightarrow GL\left( \mathcal{H}_{gL^{\circ
}}\right) $, given by $\pi _{L^{\circ }}^{g}\left( h\right) =\pi _{gL^{\circ
}}\left( g\left( h\right) \right) $. The equivariant structure is the
isomorphisms of $H\left( V\right) $-vector bundles 
\begin{equation}
\theta _{g}:g^{\ast }\mathfrak{H\rightarrow H},  \label{equivariant_eq}
\end{equation}%
which on the level of fibers, sends $f\in \mathcal{H}_{gL^{\circ }}$ to $%
f\circ g\in \mathcal{H}_{L^{\circ }}$.

\subsection{The strong Stone--von Neumann property \label{SS-vN_sub}}

By Theorem \ref{S-vN_thm}, for every pair $\left( M^{\circ },L^{\circ
}\right) \in OLag\left( V\right) ^{2}$, the models $\mathcal{H}_{L^{\circ }}$
and $\mathcal{H}_{M^{\circ }}$ are isomorphic as representations of $H\left(
V\right) $, moreover, since the Heisenberg representation is irreducible,
the vector space $\mathrm{Hom}_{H}\left( \mathcal{H}_{L^{\circ }},\mathcal{H}%
_{M}^{\circ }\right) $ of intertwining morphisms is one dimensional.
Roughly, the strong Stone--von Neumann property asserts the existence of a
distinguished element $F_{M^{\circ },L^{\circ }}\in $ $\mathrm{Hom}%
_{H}\left( \mathcal{H}_{L^{\circ }},\mathcal{H}_{M^{\circ }}\right) $, for
every pair $\left( M^{\circ },L^{\circ }\right) \in OLag\left( V\right) ^{2}$%
. The precise statement involves the following definition:

\begin{definition}
Let $\mathfrak{E}\rightarrow OLag\left( V\right) $ be an $H\left( V\right)
^{k}$-vector bundle. A \underline{trivialization} of $\mathfrak{E}$ is a
system of intertwining isomorphisms $\{E_{M^{\circ },L^{\circ }}\in \mathrm{%
Hom}_{H\left( V\right) ^{k}}(\mathfrak{E}_{L^{\circ }},\mathfrak{E}%
_{M^{\circ }}):\left( M^{\circ },L^{\circ }\right) \in OLag\left( V\right)
^{2}\}$ satisfying the following multiplicativity condition%
\begin{equation*}
E_{N^{\circ },M^{\circ }}\circ E_{M^{\circ },L^{\circ }}=E_{N^{\circ
},L^{\circ }},
\end{equation*}%
for every $N^{\circ },M^{\circ },L^{\circ }\in OLag\left( V\right) $.
\end{definition}

\begin{remark}
Intuitively, a trivialization of a $H\left( V\right) ^{k}$-vector bundle $%
\mathfrak{E}\rightarrow OLag\left( V\right) $ can be thought of as a flat
connection, compatible with the Heisenberg action and admitting a trivial
monodromy.
\end{remark}

\begin{theorem}[The strong S-vN property]
\label{SS-vN_thm}The $H\left( V\right) $-vector bundle $\mathfrak{H}$ admits
a natural trivialization $\{T_{M^{\circ },L^{\circ }}\}$.
\end{theorem}

The intertwining morphisms $T_{M^{\circ },L^{\circ }}$ in the above theorem
will be referred to as \textit{the canonical intertwining morphisms}.

\subsection{A linear algebra proof of the strong S-vN property}

The proof of Theorem \ref{SS-vN_thm} proceeds in several steps. First we
note that the vector bundle $\mathfrak{H}$ admits a natural partial
trivialization: Let us denote by $U_{2}\subset OLag\left( V\right) ^{2}$ the
subset consisting of pairs of oriented Lagrangians $\left( M^{\circ
},L^{\circ }\right) \in OLag\left( V\right) ^{2}$ which are in \underline{%
general position}, that is $L+M=V$. For every $\left( M^{\circ },L^{\circ
}\right) \in U_{2}$, define the intertwining morphisms 
\begin{equation}
T_{M^{\circ },L^{\circ }}=A_{M^{\circ },L^{\circ }}\cdot F_{M^{\circ
},L^{\circ }},  \label{ansatz_eq}
\end{equation}%
where $F_{M^{\circ },L^{\circ }}:\mathcal{H}_{L^{\circ }}\rightarrow 
\mathcal{H}_{M^{\circ }}$ is the averaging morphism given by%
\begin{equation*}
F_{M^{\circ },L^{\circ }}\left[ f\right] \left( h\right) =\sum \limits_{m\in
M}f\left( m\cdot h\right) ,
\end{equation*}%
for every $f\in \mathcal{H}_{L^{\circ }}$ and $A_{M^{\circ },L^{\circ }}$ is
a normalization constant given by 
\begin{equation}
A_{M^{\circ },L^{\circ }}=\left( G_{1}/q\right) ^{n}\sigma (\left( -1\right)
^{\left( \QATOP{n}{2}\right) }\omega _{\wedge }\left( o_{L},o_{M}\right) ),
\label{normal-ansatz_eq}
\end{equation}

where

\begin{itemize}
\item $\sigma $ is the unique quadratic character (also called the Legendre
character) of the multiplicative group $G_{m}=\mathbb{F}_{q}^{\times }$.

\item $G_{1}$ is the one dimensional Gauss sum 
\begin{equation*}
G_{1}=\sum \limits_{z\in \mathbb{F}_{q}}\psi (\tfrac{1}{2}z^{2}).
\end{equation*}

\item $\omega _{\wedge }:\bigwedge \nolimits^{top}L\times $ $\bigwedge
\nolimits^{top}M\rightarrow $ $\mathbb{F}_{q}$ is the pairing induced by the
symplectic form.
\end{itemize}

Let us denote by $U_{3}\subset OLag\left( V\right) ^{3}$ the subset
consisting of triples of oriented Lagrangians $\left( N^{\circ },M^{\circ
},L^{\circ }\right) $ which are in general position pairwisely.

\begin{proposition}
\label{multiplicativity_prop}For every $\left( N^{\circ },M^{\circ
},L^{\circ }\right) \in U_{3}$%
\begin{equation*}
T_{N^{\circ },L^{\circ }}=T_{N^{\circ },M^{\circ }}\circ T_{M^{\circ
},L^{\circ }}\text{.}
\end{equation*}
\end{proposition}

For a proof, see Appendix \ref{proofs_sec}.

Theorem \ref{SS-vN_thm}, now, follows from

\begin{proposition}
\label{extension_prop}The sub-system $\left \{ T_{M^{\circ },L^{\circ
}}:\left( M^{\circ },L^{\circ }\right) \in U_{2}\right \} $ extends in a
unique manner to a trivialization of $\mathfrak{H}$.
\end{proposition}

For a proof, see Appendix \ref{proofs_sec}.

We will refer to $\left \{ T_{M^{\circ },L^{\circ }}\right \} $ as the
system of \textit{canonical intertwining morphisms}.

\subsection{Explicit formulas for the canonical intertwining morphisms}

The canonical intertwining morphism $T_{M^{\circ },L^{\circ }}$\ can be
written in a closed form for a general pair $\left( M^{\circ },L^{\circ
}\right) \in OLag\left( V\right) ^{2}$. In order to do that we need to fix
some additional terminology: Denote $I=M\cap L$ and let $n_{I}=\frac{\dim
\left( I^{\perp }/I\right) }{2}$.

There are canonical decompositions 
\begin{eqnarray*}
\tbigwedge \nolimits^{top}M &=&\tbigwedge \nolimits^{top}I\tbigotimes
\tbigwedge \nolimits^{top}M/I, \\
\tbigwedge \nolimits^{top}L &=&\tbigwedge \nolimits^{top}I\tbigotimes
\tbigwedge \nolimits^{top}L/I.
\end{eqnarray*}

In terms of these decompositions, the orientations on $M$ and $L$ can be
written in the form $o_{M}=\iota _{M}\otimes o_{M/I}$ and $o_{L}=\iota
_{L}\otimes o_{L/I}$ respectively. Let $F_{M^{\circ },L^{\circ }}:\mathcal{H}%
_{L^{\circ }}\rightarrow \mathcal{H}_{M^{\circ }}$ denote the averaging
morphism 
\begin{equation*}
F_{M^{\circ },L^{\circ }}\left[ f\right] \left( h\right) =\sum \limits_{%
\overline{m}\in M/I}f\left( m\cdot h\right) ,
\end{equation*}%
for every $f\in \mathcal{H}_{L^{\circ }}$ where in the above summation, $%
m\in M$ is any element lying over $\overline{m}\in M/I$. Define the
normalization constant 
\begin{equation*}
A_{M^{\circ },L^{\circ }}=\left( G_{1}/q\right) ^{n_{I}}\sigma \left( \left(
-1\right) ^{\left( \QATOP{n_{I}}{2}\right) }\frac{\iota _{L}}{\iota _{M}}%
\cdot \omega _{\wedge }\left( o_{L/I},o_{M/I}\right) \right) ,
\end{equation*}%
where $\omega _{\wedge }:\bigwedge \nolimits^{top}L/I\times $ $\bigwedge
\nolimits^{top}M/I\rightarrow $ $\mathbb{F}_{q}$ is the pairing induced by
the symplectic form $\omega $.

\begin{proposition}
\label{explicit_prop} For every $\left( M^{\circ },L^{\circ }\right) \in
OLag\left( V\right) ^{2}$ 
\begin{equation*}
T_{M^{\circ },L^{\circ }}=A_{M^{\circ },L^{\circ }}\cdot F_{M^{\circ
},L^{\circ }}\text{.}
\end{equation*}
\end{proposition}

For a proof, see Appendix \ref{proofs_sec}.

\subsection{Kernel presentation of an intertwining morphism}

An explicit way to present an intertwining morphism is via a kernel
function. Fix a pair $\left( M^{\circ },L^{\circ }\right) \in OLag\left(
V\right) ^{2}$ of oriented Lagrangians and let $%
%TCIMACRO{\U{2102} }%
%BeginExpansion
\mathbb{C}
%EndExpansion
\left( M\backslash H\left( V\right) /L,\psi \right) \subset 
%TCIMACRO{\U{2102} }%
%BeginExpansion
\mathbb{C}
%EndExpansion
\left( H\left( V\right) ,\psi \right) $ be the subspace of functions $K\in 
%TCIMACRO{\U{2102} }%
%BeginExpansion
\mathbb{C}
%EndExpansion
\left( H\left( V\right) ,\psi \right) $,$\ $satisfying the equivariance
property $K\left( m\cdot h\cdot l\right) =K\left( h\right) $, for every $%
m\in M$ and $l\in L$. Given a function $K\in 
%TCIMACRO{\U{2102} }%
%BeginExpansion
\mathbb{C}
%EndExpansion
\left( M\backslash H\left( V\right) /L,\psi \right) $, we can associate to
it the intertwining morphism $I\left[ K\right] \in $ $\mathrm{Hom}_{H\left(
V\right) }(\mathcal{H}_{L^{\circ }},\mathcal{H}_{M^{\circ }})$ defined by 
\begin{equation*}
I\left[ K\right] \left( f\right) =K\ast f=m_{!}\left( K\boxtimes _{Z\cdot
L}f\right) ,
\end{equation*}%
for\ every $f\in \mathcal{H}_{L^{\circ }}$, where $K\boxtimes _{Z\cdot L}f$ $%
\ $is the descent of the function $K\boxtimes f\in 
%TCIMACRO{\U{2102} }%
%BeginExpansion
\mathbb{C}
%EndExpansion
\left( H\left( V\right) \times H\left( V\right) \right) $ to $H\left(
V\right) \times _{Z\cdot L}H\left( V\right) $---the quotient of $H\left(
V\right) \times H\left( V\right) $ by the action $x\cdot
(h_{1},h_{2})=(h_{1}x,x^{-1}h_{2})$ for $x\in Z\cdot L$---and $m_{!}$
denotes the operation of summation along the fibers of the multiplication
mapping $m:H\left( V\right) \times H\left( V\right) \rightarrow H\left(
V\right) $. We call the function $K$ an\textit{\ intertwining kernel}. The
procedure that we just described defines an isomorphism of vector spaces 
\begin{equation*}
I:%
%TCIMACRO{\U{2102} }%
%BeginExpansion
\mathbb{C}
%EndExpansion
\left( M\backslash H\left( V\right) /L,\psi \right) \longrightarrow \mathrm{%
Hom}_{H\left( V\right) }(\mathcal{H}_{L^{\circ }},\mathcal{H}_{M^{\circ }}).
\end{equation*}

Fix a triple $\left( N^{\circ },M^{\circ },L^{\circ }\right) \in OLag\left(
V\right) ^{3}$. Given kernels $K_{1}\in 
%TCIMACRO{\U{2102} }%
%BeginExpansion
\mathbb{C}
%EndExpansion
\left( N\backslash H\left( V\right) /M,\psi \right) $ and $K_{2}\in 
%TCIMACRO{\U{2102} }%
%BeginExpansion
\mathbb{C}
%EndExpansion
\left( M\backslash H\left( V\right) /L,\psi \right) $, their convolution $%
K_{1}\ast K_{2}=m_{!}\left( K_{1}\boxtimes _{Z\cdot M}K_{2}\right) $ lies in 
$%
%TCIMACRO{\U{2102} }%
%BeginExpansion
\mathbb{C}
%EndExpansion
\left( N\backslash H\left( V\right) /L,\psi \right) $. Moreover, the
transform $I$ sends convolution of kernels to composition of operators 
\begin{equation*}
I\left[ K_{1}\ast K_{2}\right] =I\left[ K_{1}\right] \circ I\left[ K_{2}%
\right] .
\end{equation*}

\subsubsection{Canonical system of intertwining kernels \label%
{sys_kernels_subsub}}

Below, we formulate a version of Theorem \ref{SS-vN_thm} in the setting of
kernels.

\begin{definition}
A system $\{K_{M^{\circ },L^{\circ }}\in 
%TCIMACRO{\U{2102} }%
%BeginExpansion
\mathbb{C}
%EndExpansion
\left( M\backslash H\left( V\right) /L,\psi \right) :\left( M^{\circ
},L^{\circ }\right) \in OLag\left( V\right) ^{2}\}$ of kernels is called 
\underline{multiplicative} if for every triple $\left( N^{\circ },M^{\circ
},L^{\circ }\right) \in OLag\left( V\right) ^{3}$ the following equation
holds 
\begin{equation*}
K_{N^{\circ },L^{\circ }}=K_{N^{\circ },M^{\circ }}\ast K_{M^{\circ
},L^{\circ }}
\end{equation*}
\end{definition}

A multiplicative system of kernels $\{K_{M^{\circ },L^{\circ }}\}$ can be
equivalently thought of as a single function $K\in 
%TCIMACRO{\U{2102} }%
%BeginExpansion
\mathbb{C}
%EndExpansion
(OLag\left( V\right) ^{2}\times H\left( V\right) )$ defined by $K\left(
M^{\circ },L^{\circ },-\right) =K_{M^{\circ },L^{\circ }}\left( -\right) $
satisfying\ the following multiplicativity relation on $OLag\left( V\right)
^{3}\times H\left( V\right) $ 
\begin{equation}
p_{12}^{\ast }K\ast p_{23}^{\ast }K=p_{13}^{\ast }K\text{,}  \label{mult_eq}
\end{equation}%
where $p_{ij}\left( \left( L_{1}^{\circ },L_{2}^{\circ },L_{3}^{\circ
}\right) ,h\right) =\left( \left( L_{i}^{\circ },L_{j}^{\circ }\right)
,h\right) $ are the projections on the $i,j$ copies of $OLag\left( V\right) $
and the left-hand side of (\ref{mult_eq}) means fiberwise convolution,
namely $p_{12}^{\ast }K\ast p_{23}^{\ast }K(L_{1}^{\circ },L_{2}^{\circ
},L_{3}^{\circ },-)=K\left( L_{1}^{\circ },L_{2}^{\circ },-\right) \ast
K\left( L_{2}^{\circ },L_{3}^{\circ },-\right) $.

For every $\left( M^{\circ },L^{\circ }\right) \in OLag\left( V\right) ^{2}$%
, there exists a unique kernel $K_{M^{\circ },L^{\circ }}$ such that $%
T_{M^{\circ },L^{\circ }}=I\left[ K_{M^{\circ },L^{\circ }}\right] $. We
will refer to $\left \{ K_{M^{\circ },L^{\circ }}\right \} $ as the system
of \textit{canonical intertwining kernels}. We will denote the corresponding
function on $OLag\left( V\right) ^{2}\times H\left( V\right) $ by $K$.

\begin{proposition}
\label{SS-vN2_prop}The function $K\in 
%TCIMACRO{\U{2102} }%
%BeginExpansion
\mathbb{C}
%EndExpansion
(OLag\left( V\right) ^{2}\times H\left( V\right) )$ satisfies%
\begin{equation*}
p_{12}^{\ast }K\ast p_{23}^{\ast }K=p_{13}^{\ast }K.
\end{equation*}
\end{proposition}

In case $\left( M^{\circ },L^{\circ }\right) \in U_{2}$, the kernel $%
K_{M^{\circ },L^{\circ }}$ is given by the following explicit formula%
\begin{equation}
K_{M^{\circ },L^{\circ }}=A_{M^{\circ },L^{\circ }}\cdot \widetilde{K}%
_{M^{\circ },L^{\circ }},  \label{ansatz2}
\end{equation}%
where $\widetilde{K}_{M^{\circ },L^{\circ }}=\left( \tau ^{-1}\right) ^{\ast
}\psi $ where $\tau =\tau _{M^{\circ },L^{\circ }}$ is the isomorphism given
by the composition $Z\hookrightarrow H\twoheadrightarrow M\backslash H\left(
V\right) /L$. The system $\{K_{M^{\circ },L^{\circ }}:\left( M^{\circ
},L^{\circ }\right) \in U_{2}\}$ yields a well defined function $%
K_{U_{2}}\in 
%TCIMACRO{\U{2102} }%
%BeginExpansion
\mathbb{C}
%EndExpansion
\left( U_{2}\times H\left( V\right) \right) $. In Section \ref{geom_sec}, we
will give an algebra geometric interpretation to the description of the
kernels $K_{M^{\circ },L^{\circ }}$ when $\left( M^{\circ },L^{\circ
}\right) \notin U_{2}$.

\subsection{The canonical vector space\label{CH_sub}}

Let us denote by $\mathsf{Symp}$ the category whose objects are finite
dimensional symplectic vector spaces over the finite field $\mathbb{F}_{q}$
and morphisms are linear isomorphisms of symplectic vector spaces. In
addition, let us denote by $\mathsf{Vect}$ the category of finite
dimensional complex vector spaces. Using Theorem \ref{SS-vN_thm}, we can
associate, in a functorial manner, a vector space $\mathcal{H}\left(
V\right) $ \ to a symplectic vector space $V\in \mathsf{Symp}$ as follows:
Define $\mathcal{H}\left( V\right) $ to be the space of "horizontal
sections" of the trivialized vector bundle $\mathfrak{H}\left( V\right) $ 
\begin{equation*}
\mathcal{H}\left( V\right) =\Gamma _{hor}\left( OLag\left( V\right) ,%
\mathfrak{H}\left( V\right) \right) ,
\end{equation*}%
where $\Gamma _{hor}\left( OLag\left( V\right) ,\mathfrak{H}\left( V\right)
\right) \subset \Gamma \left( OLag\left( V\right) ,\mathfrak{H}\left(
V\right) \right) $ is the sub-space consisting of sections $(v_{L^{\circ
}}\in \mathcal{H}_{L^{\circ }}:L^{\circ }\in OLag\left( V\right) )$
satisfying $T_{M^{\circ },L^{\circ }}\left( v_{L^{\circ }}\right)
=v_{M^{\circ }}$ for every $\left( M^{\circ },L^{\circ }\right) \in
OLag\left( V\right) ^{2}$. The vector space $\mathcal{H}\left( V\right) $ \
will be referred to as the \textit{canonical vector space} associated with $%
V $.

\begin{notation}
The definition of the vector space $\mathcal{H}\left( V\right) $ depends on
a choice of a central character $\psi $. Sometime, we will use the notation $%
\mathcal{H}\left( V,\psi \right) $ to emphasize this point.
\end{notation}

\begin{proposition}[Functoriality]
\label{functor_prop}The rule $V\mapsto \mathcal{H}\left( V\right) $ \
establishes a contravariant (quantization) functor 
\begin{equation*}
\mathcal{H}:\mathsf{Symp}\longrightarrow \mathsf{Vect}.
\end{equation*}
\end{proposition}

For a proof, see Appendix \ref{proofs_sec}.

Considering a fixed symplectic vector space $V$, we obtain as a consequence
a representation $\left( \rho _{V},Sp\left( V\right) ,\mathcal{H}\left(
V\right) \right) $, with $\rho _{V}\left( g\right) =\mathcal{H}\left(
g^{-1}\right) $, for every $g\in Sp\left( V\right) $. The representation $%
\rho _{V}$ is isomorphic to the Weil representation and will be referred to
as the \textit{canonical model} of the Weil representation.

\begin{remark}
The canonical model $\rho _{V}$ can be viewed from another perspective: \ We
begin with the total vector space $\Gamma \left( V\right) =\Gamma \left(
OLag\left( V\right) ,\mathfrak{H}\left( V\right) \right) $ and make the
following two observations. First observation, is that the symplectic group $%
Sp\left( V\right) $ acts naturally on $\Gamma \left( V\right) $, the action
is of a geometric nature---induced from the diagonal action on $OLag\left(
V\right) \times H\left( V\right) $. Second observation, is that the system $%
\left \{ T_{M^{\circ },L^{\circ }}\right \} $ defines an $Sp\left( V\right) $%
-invariant idempotent, which should be thought of as a total Fourier
transform, $T:\Gamma \left( V\right) \rightarrow \Gamma \left( V\right) $
given by 
\begin{equation*}
T\left( v_{L^{\circ }}\right) =\frac{1}{\# \left( OLag\left( V\right)
\right) }\tbigoplus \limits_{M^{\circ }\in OLag\left( V\right) }T_{M^{\circ
},L^{\circ }}\left( v_{L^{\circ }}\right) ,
\end{equation*}%
\ for every $L^{\circ }\in OLag\left( V\right) $ and $v_{L^{\circ }}\in 
\mathcal{H}_{L^{\circ }}$. The situation is summarized in the following
diagram%
\begin{equation*}
Sp\left( V\right) \circlearrowright \Gamma \left( V\right) \circlearrowleft
T.
\end{equation*}%
The canonical model is given by the image of $T$, that is, $\mathcal{H}%
\left( V\right) =T\left( \Gamma \left( V\right) \right) $. The nice thing
about this point of view is that it shows a clear distinction between
operators associated with action of the symplectic group and operators
associated with intertwining morphisms. Finally, we remark that we can,
also, consider the $Sp\left( V\right) $-invariant idempotent $T^{\bot }=Id-T$
and the associated representation $\left( \rho _{V}^{\bot },Sp\left(
V\right) ,\mathcal{H}\left( V\right) ^{\bot }\right) $, with $\mathcal{H}%
\left( V\right) ^{\bot }=T^{\bot }\left( \Gamma \left( V\right) \right) $.
The meaning of this representation is unclear.
\end{remark}

\subsection{Properties of the quantization functor}

\subsubsection{Compatibility with Cartesian products}

The category $\mathsf{Symp}$ \ admits a monoidal structure given by
Cartesian product of symplectic vector spaces. The category $\mathsf{Vect}$
admits the standard monoidal structure given by tensor product. The functor $%
\mathcal{H}$ is a monoidal functor with respect to these monoidal structures.

\begin{proposition}
\label{Cartesian_prop}For every $V_{1},V_{2}\in \mathsf{Symp}$, there exists
a natural isomorphism 
\begin{equation*}
\alpha _{V_{1}\times V_{2}}:\mathcal{H}\left( V_{1}\times V_{2}\right) 
\overset{\simeq }{\rightarrow }\mathcal{H}\left( V_{1}\right) \otimes 
\mathcal{H}\left( V_{2}\right) \text{.}
\end{equation*}
\end{proposition}

For a proof, see Appendix \ref{proofs_sec}.

As a result, we obtain the following compatibility condition between the
canonical models of the Weil representation 
\begin{equation}
\alpha _{V_{1}\times V_{2}}:\left( \rho _{V_{1}\times V_{2}}\right)
_{|Sp\left( V_{1}\right) \times Sp\left( V_{2}\right) }\overset{\simeq }{%
\longrightarrow }\rho _{V_{1}}\otimes \rho _{V_{2}}.  \label{product_eq}
\end{equation}

\begin{remark}
Condition (\ref{product_eq}) has an interesting consequence in case the
ground field is $\mathbb{F}_{3}$ \cite{G}. In this case, the group $Sp\left(
V\right) $ is not perfect when $\dim V=2$, therefore, in this particular
situation, the Weil representation is not uniquely defined. Nevertheless,
since the group $Sp\left( V\right) $ becomes perfect when $\dim V>2$, the
canonical model gives a natural choice for the Weil representation in the
"singular" dimension, $\dim V=2$.
\end{remark}

\subsubsection{Compatibility with symplectic duality}

Let $V=\left( V,\omega \right) \in \mathsf{Symp}$ and let us denote by $%
\overline{V}=\left( V,-\omega \right) $ the symplectic dual of $V$.

\begin{proposition}
\label{duality_prop}There exists a natural non-degenerate pairing 
\begin{equation*}
\left \langle \cdot ,\cdot \right \rangle _{V}:\mathcal{H}\left( \overline{V}%
\right) \times \mathcal{H}\left( V\right) \rightarrow 
%TCIMACRO{\U{2102} }%
%BeginExpansion
\mathbb{C}
%EndExpansion
,
\end{equation*}

where $V\in \mathsf{Symp}$.
\end{proposition}

For a proof, see Appendix \ref{proofs_sec}.

\subsubsection{Compatibility with symplectic reduction}

Let $V\in \mathsf{Symp}$. Let $I$ be an isotropic subspace in $V$ considered
as an abelian subgroup in $H\left( V\right) $. The fiberwise action of $%
H\left( V\right) $ on the vector bundle $\mathfrak{H}$ induces an action of $%
H\left( V\right) $ on $\mathcal{H}\left( V\right) $, using this action we
can associate to $I$ the subspace $\mathcal{H}\left( V\right) ^{I}$ of $I$%
-invariant vectors. In addition, we can form the symplectic reduction%
\footnote{%
Note that since $I$ is isotropic, $I\subset I^{\perp }$ and $I^{\perp }/I$
is equipped with a natural symplectic structure.
\par
{}} $I^{\bot }/I$ and consider the canonical vector space $\mathcal{H}\left(
I^{\bot }/I\right) $. Roughly, we claim that the vector spaces $\mathcal{H}%
\left( I^{\bot }/I\right) $ and $\mathcal{H}\left( V\right) ^{I}$ are
naturally isomorphic. The precise statement involves the following definition

\begin{definition}
An \underline{oriented isotropic subspace} in $V$ is a pair $I^{\circ
}=\left( I,o_{I}\right) $, where $I\subset V$ is an isotropic subspace and $%
o_{I}\in \tbigwedge \nolimits^{top}I$ is a non-trivial vector.
\end{definition}

\begin{proposition}
\label{reduction_prop}There exists a natural isomorphism%
\begin{equation*}
\alpha _{\left( I^{\circ },V\right) }:\mathcal{H}\left( V\right) ^{I}\overset%
{\simeq }{\rightarrow }\mathcal{H}\left( I^{\perp }/I\right) ,
\end{equation*}%
where, $V\in \mathsf{Symp}$ and $I^{\circ }$ an oriented isotropic subspace
in $V$. The naturality condition is 
\begin{equation*}
\mathcal{H}\left( f_{I}\right) \circ \alpha _{\left( J^{\circ },U\right)
}=\alpha _{\left( I^{\circ },V\right) }\circ \mathcal{H}\left( f\right) ,
\end{equation*}%
for every $f\in \mathrm{Mor}_{\mathsf{Symp}}\left( V,U\right) $ such that $%
f\left( I^{\circ }\right) =J^{\circ }$ and $f_{I}\in \mathrm{Mor}_{\mathsf{%
Symp}}\left( I^{\perp }/I,J^{\perp }/J\right) $ is the induced morphism.
\end{proposition}

For a proof, see Appendix \ref{proofs_sec}.

As a result, we obtain the following compatibility condition between the
canonical models of the Weil representation: Fix $V\in $ $\mathsf{Symp}$ and
let $I^{\circ }$ be an enhanced isotropic subspace in $V$. Let $P\subset
Sp\left( V\right) $ be the subgroup of elements $g\in Sp\left( V\right) $
such that $g\left( I^{\circ }\right) =I^{\circ }$. The isomorphism $\alpha
_{\left( I^{\circ },V\right) }$ establishes the following isomorphism 
\begin{equation*}
\alpha _{\left( I^{\circ },V\right) }:\left( \rho _{V}\right) _{|P}\overset{%
\simeq }{\longrightarrow }\rho _{I^{\perp }/I}\circ \pi ,
\end{equation*}%
where $\pi :P\rightarrow Sp\left( I^{\perp }/I\right) $ is the canonical
homomorphism.

\section{Geometric canonical intertwining kernels\label{geom_sec}}

In this section, we construct a geometric counterpart to the set-theoretic
system of canonical intertwining kernels. In particular, we obtain an
algebra-geometric interpretation for the kernels $K_{M^{\circ },L^{\circ }}$
when $\left( M^{\circ },L^{\circ }\right) \notin U_{2},$ and an alternative
proof for the strong S-vN theorem.

\subsection{Preliminaries from algebraic geometry}

We denote by $k$ an algebraic closure of the field $\mathbb{F}_{q}$. Next we
have to use some space to recall notions and notations from algebraic
geometry and the theory of $\ell $-adic sheaves. \ 

\subsubsection{Varieties}

In the sequel, we are going to translate back and forth between algebraic
varieties defined over the finite field $\mathbb{F}_{q}$ and their
corresponding sets of rational points. In order to prevent confusion between
the two, we use bold-face letters to denote a variety $\mathbf{X}$ and
normal letters $X$ to denote its corresponding set of rational points $X=%
\mathbf{X}(\mathbb{F}_{q})$. For us, a variety $\mathbf{X}$ over the finite
field is a quasi projective algebraic variety, such that the defining
equations are given by homogeneous polynomials with coefficients in the
finite field $\mathbb{F}_{q}$. In this situation, there exists a (geometric) 
\textit{Frobenius} endomorphism $Fr:\mathbf{X\rightarrow X}$, which is a
morphism of algebraic varieties. We denote by $X$ \ the set of points fixed
by $Fr$, i.e., 
\begin{equation*}
X=\mathbf{X}(\mathbb{F}_{q})=\mathbf{X}^{Fr}=\{x\in \mathbf{X}:Fr(x)=x\}.
\end{equation*}

The category of algebraic varieties over $\mathbb{F}_{q}$ will be denoted by 
$\mathsf{Var}_{\mathbb{F}_{q}}$.

\subsubsection{Sheaves}

Let $\mathsf{D}^{b}(\mathbf{X)}$ denote the bounded derived category of
constructible $\ell $-adic sheaves on $\mathbf{X}$ \cite{BBD}. We denote by $%
\mathsf{Perv}(\mathbf{X)}$ the Abelian category of perverse sheaves on the
variety $\mathbf{X}$, that is the heart with respect to the autodual
perverse t-structure in $\mathsf{D}^{b}(\mathbf{X})$. An object $\mathcal{%
F\in }\mathsf{D}^{b}(\mathbf{X)}$ is called $N$-perverse if $\mathcal{F[}%
N]\in \mathsf{Perv}(\mathbf{X)}$. Finally, we recall the notion of a Weil
structure (Frobenius structure) \cite{D}. A Weil structure associated to an
object $\mathcal{F\in }\mathsf{D}^{b}(\mathbf{X)}$ is an isomorphism%
\begin{equation*}
\theta :Fr^{\ast }\mathcal{F}\overset{\sim }{\longrightarrow }\mathcal{F}%
\text{.}
\end{equation*}

A pair $(\mathcal{F},\theta )$ is called a Weil object. By an abuse of
notation we often denote $\theta $ also by $Fr$. We choose once an
identification $\overline{%
%TCIMACRO{\U{211a} }%
%BeginExpansion
\mathbb{Q}
%EndExpansion
}_{\ell }\simeq 
%TCIMACRO{\U{2102} }%
%BeginExpansion
\mathbb{C}
%EndExpansion
$, hence all sheaves are considered over the complex numbers.

\begin{remark}
All the results in this section make perfect sense over the field $\overline{%
%TCIMACRO{\U{211a} }%
%BeginExpansion
\mathbb{Q}
%EndExpansion
}_{\ell }$, in this respect the identification of $\overline{%
%TCIMACRO{\U{211a} }%
%BeginExpansion
\mathbb{Q}
%EndExpansion
}_{\ell }$ with $%
%TCIMACRO{\U{2102} }%
%BeginExpansion
\mathbb{C}
%EndExpansion
$ \ is redundant. The reason it is specified is in order to relate our
results with the standard constructions.
\end{remark}

Given a Weil object $(\mathcal{F},Fr^{\ast }\mathcal{F\simeq F})$ one can
associate to it a function $f^{\mathcal{F}}:X\rightarrow 
%TCIMACRO{\U{2102} }%
%BeginExpansion
\mathbb{C}
%EndExpansion
$ to $\mathcal{F}$ as follows 
\begin{equation*}
f^{\mathcal{F}}(x)=\tsum \limits_{i}(-1)^{i}Tr(Fr_{|H^{i}(\mathcal{F}_{x})}).
\end{equation*}

This procedure is called \textit{Grothendieck's sheaf-to-function
correspondence\cite{Gr}}. Another common notation for the function $f^{%
\mathcal{F}}$ is $\chi _{Fr}(\mathcal{F)}$, which is called the \textit{%
Euler characteristic} of the sheaf $\mathcal{F}.$

\subsection{Geometrization}

We shall now start the geometrization procedure.

\subsubsection{Replacing sets by varieties}

The first step we take is to replace all sets involved by their geometric
counterparts, i.e., algebraic varieties. The symplectic space $(V,\omega )$
is naturally identified as the set $V=\mathbf{V}(\mathbb{F}_{q})$, where $%
\mathbf{V}$ is a $2n$-dimensional symplectic vector space in $\mathsf{Var}_{%
\mathbb{F}_{q}}$. The Heisenberg group $H$ is naturally identified as the
set $H=\mathbf{H}(\mathbb{F}_{q})$, where $\mathbf{H}=\mathbf{V\times }%
\mathbb{G}_{a}$\ is the corresponding group variety. Finally, $OLag\left(
V\right) =\mathbf{OLag}(\mathbf{V})\left( \mathbb{F}_{q}\right) $, where $%
\mathbf{OLag}\left( \mathbf{V}\right) $ is the variety of oriented
Lagrangians in $\mathbf{V}$.

\subsubsection{Replacing functions by sheaves}

The second step is to replace functions by their sheaf-theoretic
counterparts \cite{Ga}. The additive character $\psi :\mathbb{F}%
_{q}\longrightarrow 
%TCIMACRO{\U{2102} }%
%BeginExpansion
\mathbb{C}
%EndExpansion
^{\times }$ is associated via the sheaf-to-function correspondence to the
Artin-Schreier sheaf $\mathcal{L}_{\psi }$ on the variety $\mathbb{G}_{a}$,
i.e., we have $f^{\mathcal{L}_{\psi }}=\psi .$ The Legendre character $%
\sigma $ on $\mathbb{F}_{q}^{\times }\simeq $ $\mathbb{G}_{m}(\mathbb{F}%
_{q}) $ is associated to the Kummer sheaf $\mathcal{L}_{\sigma }$ on the
variety $\mathbb{G}_{m}$. The one dimensional Gauss sum $G_{1}$ is
associated with the Weil object 
\begin{equation*}
\mathcal{G}_{1}\mathcal{=}\tint \limits_{\mathbb{G}_{a}}\mathcal{L}_{\psi (%
\frac{1}{2}z^{2})}\in \mathsf{D}^{b}(\mathbf{pt}),
\end{equation*}%
where, for the rest of this paper, $\int =\int_{!}$ denotes integration with
compact support \cite{BBD}. Grothendieck's Lefschetz trace formula \cite{Gr}
implies that, indeed, $f^{\mathcal{G}_{1}}=G_{1}.$ In fact, there exists a
quasi-isomorphism $\mathcal{G}_{1}\overset{q.i}{\longrightarrow }$ $H^{1}(%
\mathcal{G}_{1}\mathcal{)}[-1]$ and $\dim H^{1}(\mathcal{G}_{1}\mathcal{)=}1$%
, hence, $\mathcal{G}_{1}$ can be thought of as a one-dimensional vector
space, equipped with a Frobenius operator, sitting at cohomological degree $%
1.$

\subsection{Canonical system of geometric intertwining kernels}

Let $\mathbf{U}_{2}\mathbf{\subset OLag}\left( \mathbf{V}\right) ^{2}$ be
the open subvariety consisting of pairs $\left( M^{\circ },L^{\circ }\right)
\in \mathbf{OLag}\left( \mathbf{V}\right) ^{2}$ which are in general
position. We define a sheaf \ "of kernels" \ $\mathcal{K}_{\mathbf{U}_{2}}$
on the variety $\mathbf{U}_{2}\times \mathbf{H}\left( \mathbf{V}\right) $ as
follows:

\begin{equation*}
\mathcal{K}_{\mathbf{U}_{2}}=\mathcal{A\otimes }\widetilde{\mathcal{K}}_{%
\mathbf{U}_{2}}.
\end{equation*}

where

\begin{itemize}
\item $\widetilde{\mathcal{K}}_{\mathbf{U}_{2}}$ is the sheaf of
non-normalized kernels\textit{\ }given by\textit{\ }%
\begin{equation*}
\widetilde{\mathcal{K}}_{\mathbf{U}_{2}}\left( M^{\circ },L^{\circ }\right)
=\left( \mathcal{\tau }^{-1}\right) ^{\ast }\mathcal{L}_{\psi },
\end{equation*}%
where $\tau =\tau _{M^{\circ },L^{\circ }}$ is the isomorphism given by the
composition 
\begin{equation*}
\mathbf{Z\hookrightarrow H\twoheadrightarrow M\backslash H/L}.
\end{equation*}

\item $\mathcal{A}$ is the "Normalization coefficient" sheaf given by%
\begin{equation}
\mathcal{A}\left( M^{\circ },L^{\circ }\right) \mathcal{=G}_{1}^{\otimes
n}\otimes \mathcal{L}_{\sigma }\left( \left( -1\right) ^{\left( \QATOP{n}{2}%
\right) }\omega _{\wedge }\left( o_{L},o_{M}\right) \right) [2n]\left(
n\right) .  \label{norml_coef}
\end{equation}
\end{itemize}

Let $n_{k}=\dim (\mathbf{OLag}\left( \mathbf{V}\right) ^{k})+n+1$ for $k\in 
%TCIMACRO{\U{2115} }%
%BeginExpansion
\mathbb{N}
%EndExpansion
$. By construction, the sheaf $\mathcal{K}_{\mathbf{U}_{2}}$ is a local
system (liss\'{e}) of rank $1$, normalized to sit in cohomological degree $%
-n_{2}$ and consequently of pure weight $w(\mathcal{K}_{\mathbf{U}_{2}})=0.$

\begin{proposition}
\label{extension3_prop}The local system $\mathcal{K}_{\mathbf{U}_{2}}$ can
be extended in a unique manner to a geometrically irreducible $[n_{2}\mathbf{%
]}$-perverse Weil sheaf $\mathcal{K}$ on $\mathbf{OLag}\left( \mathbf{V}%
\right) ^{2}\mathbf{\times H}\left( \mathbf{V}\right) $ of pure weight $w(%
\mathcal{K})=0$. Moreover, there exists an isomorphism 
\begin{equation*}
p_{13}^{\ast }\mathcal{K\simeq }p_{12}^{\ast }\mathcal{K\ast }p_{23}^{\ast }%
\mathcal{K}.
\end{equation*}
\end{proposition}

The proof of Proposition \ref{extension3_prop} proceeds in several steps.
The construction of $\mathcal{K}$ uses the functor of middle extension,
namely, take 
\begin{equation}
\mathcal{K=}j_{!\ast }\mathcal{K}_{\mathbf{U}_{2}},  \label{kernelsheaf}
\end{equation}%
where $j:\mathbf{U}_{2}\mathbf{\times H}\left( \mathbf{V}\right) \mathbf{%
\hookrightarrow OLag}\left( \mathbf{V}\right) ^{2}\times \mathbf{H}\left( 
\mathbf{V}\right) $ is the open imbedding and\ $j_{!\ast \text{ }}$is the
functor of middle extension \cite{BBD}. It follows directly from the
construction that the sheaf $\  \mathcal{K}$ is irreducible $[n_{2}]$%
-perverse of pure weight $0$.

\begin{remark}
Proposition \ref{extension3_prop} establishes an alternative proof of
Proposition \ref{SS-vN2_prop} as follows: Let $K=f^{\mathcal{K}}$. The
multiplicativity condition for the sheaf $\mathcal{K}$ implies that $K$ is
multiplicative. Moreover, since $K_{U_{2}}=f^{\mathcal{K}_{\mathbf{U}_{2}}}$%
, it implies that $K$ extends the function $K_{U_{2}}$. In addition, the
geometric construction yields a description of the kernels $K_{M^{\circ
},L^{\circ }}$ when $\left( M^{\circ },L^{\circ }\right) \notin U_{2}$ in
terms of the middle extension of $\mathcal{K}_{\mathbf{U}_{2}}$.
\end{remark}

\subsection{Proof of the multiplicativity property}

Denote by $\mathbf{U}_{3}\subset \mathbf{OLag}\left( \mathbf{V}\right) ^{3}$
the open subvariety consisting of triples $\left( N^{\circ },M^{\circ
},L^{\circ }\right) $ which are in general position pairwisely.\ 

\begin{lemma}
\label{isomorphism1_lemma}There exists an isomorphism on $\mathbf{U}%
_{3}\times \mathbf{H}\left( \mathbf{V}\right) $ 
\begin{equation*}
p_{13}^{\ast }\mathcal{K\simeq }p_{12}^{\ast }\mathcal{K\ast }p_{23}^{\ast }%
\mathcal{K}.
\end{equation*}
\end{lemma}

For a proof, see Appendix \ref{proofs_sec}.

Let $\mathbf{V}_{3}\subset \mathbf{OLag}\left( \mathbf{V}\right) ^{3}$
denote the open subvariety consisting of triples $\left( N^{\circ },M^{\circ
},L^{\circ }\right) $ such that $\left( N^{\circ },M^{\circ }\right) ,\left(
M^{\circ },L^{\circ }\right) \in \mathbf{U}_{2}$.

\begin{lemma}
\label{isomorphism2_lemma}There exists an isomorphism on $\mathbf{V}%
_{3}\times \mathbf{H}\left( \mathbf{V}\right) $ 
\begin{equation*}
p_{13}^{\ast }\mathcal{K\simeq }p_{12}^{\ast }\mathcal{K\ast }p_{23}^{\ast }%
\mathcal{K}.
\end{equation*}
\end{lemma}

For a proof, see Appendix \ref{proofs_sec}.

Using Lemma \ref{isomorphism1_lemma} we conclude that the sheaves $%
p_{13}^{\ast }\mathcal{K}$ and $p_{12}^{\ast }\mathcal{K\ast }p_{23}^{\ast }%
\mathcal{K}$ are isomorphic on the open subvariety $\mathbf{U}_{3}\times 
\mathbf{H}\left( \mathbf{V}\right) $. The sheaf $p_{13}^{\ast }\mathcal{K}$
is irreducible $\left[ n_{3}\right] $-perverse as a pull-back by a smooth,
surjective with connected fibers morphism, of an irreducible $[n_{2}]$%
-perverse sheaf on $\mathbf{OLag}\left( \mathbf{V}\right) ^{2}\times \mathbf{%
H}\left( \mathbf{V}\right) $, hence, it is enough to show that the sheaf $%
p_{12}^{\ast }\mathcal{K\ast }p_{23}^{\ast }\mathcal{K}$ is irreducible $%
\left[ n_{3}\right] $-perverse.

Let us denote by $\mathbf{V}_{4}\subset \mathbf{OLag}\left( \mathbf{V}%
\right) ^{4}$ the open subvariety consisting of quadruples $\left( N^{\circ
},S^{\circ },M^{\circ },L^{\circ }\right) \in \mathbf{OLag}\left( \mathbf{V}%
\right) ^{4}$ such that $\left( N^{\circ },S^{\circ }\right) ,\left(
S^{\circ },M^{\circ }\right) \in \mathbf{U}_{2}$. The projection $p_{134}:%
\mathbf{V}_{4}\times \mathbf{H}\left( \mathbf{V}\right) \rightarrow \mathbf{%
OLag}\left( \mathbf{V}\right) ^{3}\times \mathbf{H}\left( \mathbf{V}\right) $
is smooth, surjective and admits connected fibers, therefore, it is enough
to show that the pull-back $p_{134}^{\ast }\left( p_{12}^{\ast }\mathcal{%
K\ast }p_{23}^{\ast }\mathcal{K}\right) $ is irreducible $\left[ n_{4}\right]
$-perverse. Using Lemma \ref{isomorphism2_lemma} and also invoking some
direct diagram chasing we obtains%
\begin{equation}
p_{134}^{\ast }\left( p_{12}^{\ast }\mathcal{K\ast }p_{23}^{\ast }\mathcal{K}%
\right) \mathcal{\simeq }p_{12}^{\ast }\mathcal{K}\ast p_{23}^{\ast }%
\mathcal{K}\ast p_{34}^{\ast }\mathcal{K}.  \label{formula1}
\end{equation}

The right-hand side of the above formula is principally a subsequent
application of a properly normalized, Fourier transforms (see Formula (\ref%
{Fo}) below) on $p_{34}^{\ast }\mathcal{K}$, hence by the Katz--Laumon
theorem \cite{KL} it is irreducible $\left[ n_{4}\right] $-perverse.

Let us summarize. We showed that both sheaves $p_{13}^{\ast }\mathcal{K}$
and $p_{12}^{\ast }\mathcal{K\ast }p_{23}^{\ast }\mathcal{K}$ are
irreducible $[n_{3}]$-perverse and are isomorphic on an open subvariety.
This implies that they must be isomorphic.

This concludes the proof of the multiplicativity property.

\appendix

\section{Proofs of statements \label{proofs_sec}}

\subsection{Proof of Proposition \protect \ref{multiplicativity_prop}}

For every $\left( N^{\circ },M^{\circ },L^{\circ }\right) \in U_{3}$ the
intertwining morphisms $F_{N^{\circ },M^{\circ }}\circ F_{M^{\circ
},L^{\circ }},F_{N^{\circ },L^{\circ }}\in \mathrm{Hom}_{H\left( V\right)
}\left( \mathcal{H}_{L^{\circ }},\mathcal{H}_{N^{\circ }}\right) $ are
proportional, namely%
\begin{equation*}
F_{N^{\circ },M^{\circ }}\circ F_{M^{\circ },L^{\circ }}=C\left( N^{\circ
},M^{\circ },L^{\circ }\right) \cdot F_{N^{\circ },L^{\circ }},
\end{equation*}%
for some $C\left( N^{\circ },M^{\circ },L^{\circ }\right) \in 
%TCIMACRO{\U{2102} }%
%BeginExpansion
\mathbb{C}
%EndExpansion
$. This follows from the fact that $\mathcal{H}_{L^{\circ }}$ and $\mathcal{H%
}_{N^{\circ }}$ are irreducible and isomorphic as representations of the
Heisenberg group $H\left( V\right) $.

The proof will proceed in two steps.

\textbf{Step 1. } We calculate the proportionality coefficient $C\left(
N^{\circ },M^{\circ },L^{\circ }\right) $.

Let $\delta _{L}\in \mathcal{H}_{L^{\circ }}$ denote the unique function
supported on $L\cdot Z$ normalized such that $\delta _{L}\left( 0\right) =1$%
. On the one hand, it is easy to verify that $F_{N^{\circ },L^{\circ
}}\left( \delta _{L}\right) \left( 0\right) =1$ therefore 
\begin{equation*}
C\left( N^{\circ },M^{\circ },L^{\circ }\right) =F_{N^{\circ },M^{\circ
}}\circ F_{M^{\circ },L^{\circ }}\left( \delta _{L}\right) \left( 0\right) 
\text{.}
\end{equation*}

On the other hand, explicit calculation reveals that 
\begin{equation*}
F_{N^{\circ },M^{\circ }}\circ F_{M^{\circ },L^{\circ }}\left( \delta
_{L}\right) \left( 0\right) =\tsum \limits_{m\in M}\psi (\tfrac{1}{2}\omega
\left( r^{L}\left( m\right) ,m\right) ),
\end{equation*}%
where $r^{L}:M\rightarrow N$ is the linear map defined by the condition $%
r^{L}\left( m\right) -m\in L$ for every $m\in M$. Furthermore, an easy
diagonalization argument implies that 
\begin{equation*}
\tsum \limits_{m\in M}\psi (\tfrac{1}{2}\omega \left( r^{L}\left( m\right)
,m\right) )=G_{1}^{n}\cdot \sigma \left( d\left[ \omega \left( r^{L}\left(
-\right) ,-\right) \right] \right) ,
\end{equation*}%
where $d\left[ \omega \left( r^{L}\left( -\right) ,-\right) \right] \in 
\mathbb{F}_{q}^{\times }/\mathbb{F}_{q}^{\times 2}$ denotes the discriminant
of the symmetric form $\omega \left( r^{L}\left( -\right) ,-\right) $. The
map $r^{L}:M\rightarrow N$ induces a map $r_{\wedge }^{L}:\tbigwedge
\nolimits^{top}M\rightarrow \tbigwedge \nolimits^{top}N$.

\begin{lemma}
\label{discr_lemma}We have%
\begin{equation*}
d\left[ \omega \left( r^{L}\left( -\right) ,-\right) \right] =\left(
-1\right) ^{\left( \QATOP{n}{2}\right) }\omega _{\wedge }\left( r_{\wedge
}^{L}\left( o_{M}\right) ,o_{M}\right) \text{.}
\end{equation*}
\end{lemma}

Summarizing, we get that 
\begin{equation}
C\left( N^{\circ },M^{\circ },L^{\circ }\right) =G_{1}^{n}\cdot \sigma
(\left( -1\right) ^{\left( \QATOP{n}{2}\right) }\omega _{\wedge }\left(
r_{\wedge }^{L}\left( o_{M}\right) ,o_{M}\right) )\text{.}
\label{cocycle_eq}
\end{equation}

\textbf{Step 2. }Denote $A\left( N^{\circ },M^{\circ },L^{\circ }\right)
=A_{N^{\circ },M^{\circ }}\cdot A_{M^{\circ },L^{\circ }}\cdot A_{N^{\circ
},L^{\circ }}^{-1}$. We will show that 
\begin{equation*}
A\left( N^{\circ },M^{\circ },L^{\circ }\right) =C\left( N^{\circ },M^{\circ
},L^{\circ }\right) ^{-1}\text{.}
\end{equation*}

Using Formula (\ref{normal-ansatz_eq}) for the normalization coefficients,
we can write $A\left( N^{\circ },M^{\circ },L^{\circ }\right) $ in the form%
\begin{equation*}
\left( G_{1}/q\right) ^{n}\sigma (\left( -1\right) ^{\left( \QATOP{n}{2}%
\right) }\omega _{\wedge }\left( o_{M},o_{N}\right) \omega _{\wedge }\left(
o_{L},o_{M}\right) \omega _{\wedge }\left( o_{L},o_{N}\right) ^{-1}).
\end{equation*}

\begin{lemma}
\label{identity_lemma}We have%
\begin{equation*}
\omega _{\wedge }\left( r_{\wedge }^{L}\left( o_{M}\right) ,o_{M}\right)
=\left( -1\right) ^{n}\omega _{\wedge }\left( o_{M},o_{N}\right) \omega
_{\wedge }\left( o_{L},o_{M}\right) \omega _{\wedge }\left(
o_{L},o_{N}\right) ^{-1}.
\end{equation*}
\end{lemma}

Using Lemma \ref{identity_lemma}, we get that $A\left( N^{\circ },M^{\circ
},L^{\circ }\right) $ is equal to 
\begin{equation*}
\left( G_{1}/q\right) ^{n}\sigma (\left( -1\right) ^{\left( \QATOP{n}{2}%
\right) +n}\omega _{\wedge }\left( r_{\wedge }^{L}\left( o_{M}\right)
,o_{M}\right) )=C\left( N^{\circ },M^{\circ },L^{\circ }\right) ^{-1},
\end{equation*}

where in the second equality we use the identity $G_{1}^{2n}=q^{n}\sigma
\left( \left( -1\right) ^{n}\right) $.

This concludes the proof of the proposition.

\subsubsection{Proof of Lemma \protect \ref{discr_lemma}}

Choose an isomorphism $\varphi :M\overset{\simeq }{\rightarrow }\mathbb{F}%
_{q}^{n}$. Let $B:M\times M\rightarrow \mathbb{F}_{q}$ be the symmetric form
defined by $B=\varphi ^{\ast }\left( \tsum x_{i}y_{i}\right) $. Present the
symmetric form $\omega \left( r^{L}\left( -\right) ,-\right) $ as 
\begin{equation*}
\omega \left( r^{L}\left( -\right) ,-\right) =B\left( A\left( -\right)
,-\right) ,
\end{equation*}%
where $A\in Mat_{n\times n}\left( \mathbb{F}_{q}\right) $ is a symmetric
matrix. By definition $d\left[ \omega \left( r^{L}\left( -\right) ,-\right) %
\right] =\det A$ $\func{mod}$ $\mathbb{F}_{q}^{\times 2}$.

Let $e_{1},e_{2},..,e_{n}$ be an orthonormal basis with respect to the form $%
B$. The argument now follows from 
\begin{eqnarray*}
\det A &=&\tsum \limits_{\sigma \in \Sigma _{n}}\left( -1\right) ^{\sigma
}\tprod \limits_{i=1}^{n}\omega (r^{L}\left( e_{i}\right) ,e_{\sigma \left(
i\right) }) \\
&=&\left( -1\right) ^{\left( \QATOP{n}{2}\right) }\omega _{\wedge
}(r_{\wedge }^{L}\left( e_{1}\wedge ..\wedge e_{n}\right) ,e_{1}\wedge
..\wedge e_{n}) \\
&=&\left( -1\right) ^{\left( \QATOP{n}{2}\right) }\omega _{\wedge
}(r^{L}\left( o_{M}\right) ,o_{M})\text{ }\func{mod}\mathbb{F}_{q}^{\times 2}%
\text{.}
\end{eqnarray*}

This concludes the proof of the lemma.

\subsubsection{ Proof of Lemma \protect \ref{identity_lemma}}

We can write $r_{\wedge }^{L}\left( o_{M}\right) =a\cdot o_{N}$, for some $%
a\in \mathbb{F}_{q}$. The constant $a$ can be computed as follows: On the
one hand, since $r^{L}\left( m\right) -m\in L$ , we have that $\omega
_{\wedge }(r_{\wedge }^{L}\left( o_{M}\right) ,o_{L})=\omega _{\wedge
}\left( o_{M},o_{L}\right) $. On the other hand, we have that $\omega
_{\wedge }\left( r_{\wedge }^{L}\left( o_{M}\right) ,o_{L}\right) =a\omega
_{\wedge }\left( o_{N},o_{L}\right) $, hence 
\begin{equation*}
a=\frac{\omega _{\wedge }\left( o_{M},o_{L}\right) }{\omega _{\wedge }\left(
o_{N},o_{L}\right) }.
\end{equation*}

Therefore we obtain that 
\begin{eqnarray*}
\omega _{\wedge }\left( r_{\wedge }^{L}\left( o_{M}\right) ,o_{M}\right)
&=&\omega _{\wedge }\left( o_{M},o_{L}\right) \omega _{\wedge }\left(
o_{N},o_{L}\right) ^{-1}\omega _{\wedge }\left( o_{N},o_{M}\right) \\
&=&\left( -1\right) ^{n}\omega _{\wedge }\left( o_{M},o_{N}\right) \omega
_{\wedge }\left( o_{L},o_{M}\right) \omega _{\wedge }\left(
o_{L},o_{N}\right) ^{-1},
\end{eqnarray*}%
where in the second equality we used the relation $\omega _{\wedge }\left(
o,o^{\prime }\right) =\left( -1\right) ^{n}\omega _{\wedge }\left( o^{\prime
},o\right) $. This concludes the proof of the lemma.

\subsection{Proof of Proposition \protect \ref{extension_prop}}

The trivialization of the vector bundle $\mathfrak{H}$ is constructed as
follows: Let $\left( N^{\circ },L^{\circ }\right) \in OLag\left( V\right)
^{2}$. In order to construct $T_{N^{\circ },L^{\circ }}$, choose a third $%
M^{\circ }\in OLag\left( V\right) $ such that $\left( N^{\circ },M^{\circ
}\right) ,\left( M^{\circ },L^{\circ }\right) \in U_{2}$ (such a choice
always exists). Define 
\begin{equation*}
T_{N^{\circ },L^{\circ }}=T_{N^{\circ },M^{\circ }}\circ T_{M^{\circ
},L^{\circ }}\text{,}
\end{equation*}%
where we note that both operators in the composition on the left hand side
are well defined. We are left to show that the definition of $T_{N^{\circ
},L^{\circ }}$ does not depends on the choice of $M^{\circ }$. Let $%
M_{i}^{\circ }\in OLag\left( V\right) $, $i=1,2$, such that $\left( N^{\circ
},M_{i}^{\circ }\right) ,\left( M_{i}^{\circ },L^{\circ }\right) \in U_{2}$.
We want to show that 
\begin{equation*}
T_{N^{\circ },M_{1}^{\circ }}\circ T_{M_{1}^{\circ },L^{\circ }}=T_{N^{\circ
},M_{2}^{\circ }}\circ T_{M_{2}^{\circ },L^{\circ }}.
\end{equation*}

Choose $M_{3}^{\circ }\in OLag\left( V\right) $ such that $\left(
M_{3}^{\circ },M_{i}^{\circ }\right) \in U_{2}$ for $i=1,2$ and, in
addition, $\left( M_{3}^{\circ },L^{\circ }\right) ,\left( M_{3}^{\circ
},N^{\circ }\right) \in U_{2}$. We have%
\begin{eqnarray*}
T_{N^{\circ },M_{1}^{\circ }}\circ T_{M_{1}^{\circ },L^{\circ }}
&=&T_{N^{\circ },M_{1}^{\circ }}\circ T_{M_{1}^{\circ },M_{3}^{\circ }}\circ
T_{M_{3}^{\circ },L^{\circ }} \\
&=&T_{N^{\circ },M_{3}^{\circ }}\circ T_{M_{3}^{\circ },L^{\circ }},
\end{eqnarray*}%
where the first and second equalities are the multiplicativity property for
triples which are in general position pairwisely (Proposition \ref%
{multiplicativity_prop}). In the same fashion, we show that $T_{N^{\circ
},M_{2}^{\circ }}\circ T_{M_{2}^{\circ },L^{\circ }}=T_{N^{\circ
},M_{3}^{\circ }}\circ T_{M_{3}^{\circ },L^{\circ }}$.

Finally, the same kind of argument shows that the complete system $%
\{T_{M^{\circ },L^{\circ }}:\left( M^{\circ },L^{\circ }\right) \in
OLag\left( V\right) ^{2}\}$ forms a trivialization.

This concludes the proof of the proposition.

\subsection{Proof of Proposition \protect \ref{explicit_prop}}

Denote $I=M\cap L$. \ Using the canonical decompositions $\tbigwedge
\nolimits^{top}M=\tbigwedge \nolimits^{top}I\otimes \tbigwedge
\nolimits^{top}M/I$ and $\tbigwedge \nolimits^{top}L=\tbigwedge
\nolimits^{top}I\otimes \tbigwedge \nolimits^{top}L/I$, we can write the
orientations on $M$ and $L$ in the form $o_{M}=\iota _{M}\otimes o_{M/I}$
and $o_{L}=\iota _{L}\otimes o_{L/I}$ respectively.

Choose a third $S^{\circ }\in OLag\left( V\right) $ such that $\left(
M^{\circ },S^{\circ }\right) ,\left( S^{\circ },L^{\circ }\right) \in U_{2}$%
. By Theorem \ref{extension_prop}, we conclude that $T_{M^{\circ },L^{\circ
}}=T_{M^{\circ },S^{\circ }}\circ T_{S^{\circ },L^{\circ }}$, furthermore,
using the explicit formula (\ref{ansatz_eq}), we can write the composition $%
T_{M^{\circ },S^{\circ }}\circ T_{S^{\circ },L^{\circ }}$ explicitly in the
form 
\begin{equation}
\left( G_{1}/q\right) ^{2n}\sigma \left( \omega _{\wedge }\left(
o_{S},o_{M}\right) \omega _{\wedge }\left( o_{L},o_{S}\right) \right)
F_{M^{\circ },S^{\circ }}\circ F_{S^{\circ },L^{\circ }}.
\label{composition1_eq}
\end{equation}

Direct computation reveals that%
\begin{eqnarray}
F_{M^{\circ },S^{\circ }}\circ F_{S^{\circ },L^{\circ }} &=&\left[ \tsum
\limits_{m\in M}\psi (\tfrac{1}{2}\omega (m,r^{L}(m)))\right] \cdot
F_{M^{\circ },L^{\circ }}  \label{composition2_eq} \\
&=&\left[ \#I\cdot \tsum \limits_{\overline{m}\in M/I}\psi (\tfrac{1}{2}%
\omega \left( \overline{m},r^{L}\left( \overline{m}\right) \right) )\right]
\cdot F_{M^{\circ },L^{\circ }},  \notag
\end{eqnarray}%
where $r^{L}:M\rightarrow S$ is the linear map defined by the condition $%
r^{L}\left( m\right) -m\in L$, for every $m\in M$; moreover, $ker\left(
r^{L}\right) =I$, hence, $r^{L}$ factors to an injective map $%
r^{L}:M/I\rightarrow S$. This also explains the second equality in (\ref%
{composition2_eq}).

Let us denote by $B:M/I\times M/I\rightarrow \mathbb{F}_{q}$ the symmetric
form on $M/I$ given by $B\left( \overline{m}_{1},\overline{m}_{2}\right) =%
\tfrac{1}{2}\omega \left( \overline{m}_{1},r^{L}\left( \overline{m}%
_{2}\right) \right) $. An easy diagonalization argument implies that 
\begin{equation*}
\tsum \limits_{\overline{m}\in M/I}\psi (\tfrac{1}{2}\omega \left( \overline{%
m},r^{L}\left( \overline{m}\right) \right) )=G_{1}^{n_{I}}\cdot \sigma
\left( d\left[ B\right] \right) ,
\end{equation*}%
where $d\left[ B\right] \in \mathbb{F}_{q}^{\times }/\mathbb{F}_{q}^{\times
2}$ denotes the discriminant of the symmetric form $B$.

Let us denote by $B_{\wedge }$ the form on $\tbigwedge \nolimits^{top}M/I$
induce from $B$.

\begin{lemma}
\label{discr2_lemma}We have%
\begin{equation*}
d\left[ B\right] =\left( -1\right) ^{\left( \QATOP{n_{I}}{2}\right)
}B_{\wedge }\left( o_{M/I},o_{M/I}\right) .
\end{equation*}
\end{lemma}

Summarizing, we obtain that $T_{M^{\circ },S^{\circ }}\circ T_{S^{\circ
},L^{\circ }}$ is equal to 
\begin{equation}
\left( G_{1}/q\right) ^{n_{I}}\sigma (\left( -1\right) ^{\left( \QATOP{n_{I}%
}{2}\right) +n}B_{\wedge }\left( o_{M/I},o_{M/I}\right) \omega _{\wedge
}\left( o_{S},o_{M}\right) \omega _{\wedge }\left( o_{L},o_{S}\right)
)F_{M^{\circ },L^{\circ }}\text{.}  \label{composition3_eq}
\end{equation}

Finally use

\begin{lemma}
\label{identity2_lemma}We have%
\begin{equation*}
B_{\wedge }\left( o_{M/I},o_{M/I}\right) =\left( -1\right) ^{n}\omega
_{\wedge }\left( o_{L/I},o_{M/I}\right) \omega _{\wedge }\left(
o_{S},o_{M}\right) \omega \left( o_{L},o_{S}\right) ^{-1}\frac{\iota _{L}}{%
\iota _{M}}\text{.}
\end{equation*}
\end{lemma}

in order to conclude that $T_{M^{\circ },S^{\circ }}\circ T_{S^{\circ
},L^{\circ }}$ is equal to 
\begin{equation*}
\left( G_{1}/q\right) ^{n_{I}}\sigma (\left( -1\right) ^{\left( \QATOP{n_{I}%
}{2}\right) }\omega _{\wedge }\left( o_{L/I},o_{M/I}\right) \frac{\iota _{L}%
}{\iota _{M}})F_{M^{\circ },L^{\circ }}\text{.}
\end{equation*}

This concludes the proof of the proposition.

\subsubsection{Proof of Lemma \protect \ref{discr2_lemma}}

Choose an isomorphism $\varphi :M/I\overset{\simeq }{\rightarrow }\mathbb{F}%
_{q}^{n_{I}}$. Let $B_{0}:M/I\times M/I\rightarrow \mathbb{F}_{q}$ be the
symmetric form defined by $B_{0}=\varphi ^{\ast }\left( \tsum
x_{i}y_{i}\right) $. Present the symmetric form $B$ as 
\begin{equation*}
B\left( -,-\right) =B_{0}\left( A\left( -\right) ,-\right) ,
\end{equation*}%
where $A\in Mat_{n_{I}\times n_{I}}\left( \mathbb{F}_{q}\right) $ is a
symmetric matrix. By definition $d\left[ B\right] =\det A$ $\func{mod}$ $%
\mathbb{F}_{q}^{\times 2}.$

Let $e_{1},e_{2},..,e_{n_{I}}$ be an orthonormal basis with respect to the
form $B_{0}$. The argument follows from the following equalities: 
\begin{eqnarray*}
\det A &=&\tsum \limits_{\sigma \in \Sigma _{n_{I}}}\left( -1\right)
^{\sigma }\tprod \limits_{i=1}^{n_{I}}\omega (e_{i},r^{L}\left( e_{\sigma
\left( i\right) }\right) ) \\
&=&\left( -1\right) ^{\frac{n_{I}\left( n_{I}-1\right) }{2}}\omega _{\wedge
}(e_{1}\wedge ..\wedge e_{n_{I}},r_{\wedge }^{L}\left( e_{1}\wedge ..\wedge
e_{n_{I}}\right) ) \\
&=&\left( -1\right) ^{\frac{n_{I}\left( n_{I}-1\right) }{2}}\omega _{\wedge
}(o_{M/I},r^{L}\left( o_{M/I}\right) )\text{ }\func{mod}\mathbb{F}%
_{q}^{\times 2}\text{.}
\end{eqnarray*}

This concludes the proof of the lemma.

\subsubsection{Proof of Lemma \protect \ref{identity2_lemma}}

Let us denote by $r_{\wedge }^{L}:\tbigwedge \nolimits^{top}M/I$ $%
\rightarrow \tbigwedge \nolimits^{top}r^{L}\left( M/I\right) $ the map
induced from $r^{L}:M/I\rightarrow r^{L}\left( M/I\right) $. There is a
canonical decomposition $\tbigwedge \nolimits^{top}S=\tbigwedge
\nolimits^{top}r^{L}\left( M/I\right) \otimes \tbigwedge
\nolimits^{top}S/r^{L}\left( M/I\right) $ and consequently the orientation $%
o_{S}$ can be written in the form $o_{S}=r_{\wedge }^{L}\left(
o_{M/I}\right) \otimes o$. On the one hand%
\begin{eqnarray*}
\omega _{\wedge }\left( o_{M},o_{S}\right) &=&\omega _{\wedge }\left( \iota
_{M}\otimes o_{M/I},r_{\wedge }^{L}\left( o_{M/I}\right) \otimes o\right) \\
&=&\omega _{\wedge }\left( \iota _{M},o\right) \cdot \omega _{\wedge }\left(
o_{M/I},r_{\wedge }^{L}\left( o_{M/I}\right) \right) .
\end{eqnarray*}

On the other hand

\begin{eqnarray*}
\omega _{\wedge }\left( o_{L},o_{S}\right) &=&\omega _{\wedge }\left( \iota
_{L}\otimes o_{L/I},r_{\wedge }^{L}\left( o_{M/I}\right) \otimes o\right) \\
&=&\omega _{\wedge }\left( \iota _{L},o\right) \cdot \omega _{\wedge }\left(
o_{L/I},r_{\wedge }^{L}\left( o_{M/I}\right) \right) \\
&=&\frac{\iota _{L}}{\iota _{M}}\omega _{\wedge }\left( \iota _{M},o\right)
\cdot \omega _{\wedge }\left( o_{L/I},o_{M/I}\right) .
\end{eqnarray*}

which implies that 
\begin{equation*}
\omega _{\wedge }\left( \iota _{M},o\right) =\omega _{\wedge }\left(
o_{L},o_{S}\right) \omega _{\wedge }\left( o_{L/I},o_{M/I}\right) ^{-1}\frac{%
\iota _{M}}{\iota _{L}}.
\end{equation*}

Altogether, we obtain that%
\begin{equation*}
\omega _{\wedge }\left( o_{M/I},r_{\wedge }^{L}\left( o_{M/I}\right) \right)
=\left( -1\right) ^{N}\omega _{\wedge }\left( o_{L/I},o_{M/I}\right) \omega
_{\wedge }\left( o_{S},o_{M}\right) \omega _{\wedge }\left(
o_{L},o_{S}\right) ^{-1}\frac{\iota _{L}}{\iota _{M}}.
\end{equation*}

This concludes the proof of the lemma.

\subsection{Proof of Proposition \protect \ref{functor_prop}}

Let $f\in \mathrm{Mor}_{\mathsf{Symp}}\left( V_{1},V_{2}\right) $ be an
isomorphism of symplectic vector spaces. The map $f$ \ induces a pair
isomorphisms 
\begin{eqnarray*}
r &:&H\left( V_{1}\right) \overset{\simeq }{\longrightarrow }H\left(
V_{2}\right) , \\
s &:&OLag\left( V_{1}\right) \overset{\simeq }{\longrightarrow }OLag\left(
V_{2}\right) .
\end{eqnarray*}

The first is the isomorphism of groups given by $r\left( v_{1},z\right)
=\left( f\left( v_{1}\right) ,z\right) $ and the second is the evident
induced bijection of sets. For every $L^{\circ }\in OLag\left( V_{1}\right) $%
, we have a pull-back isomorphism $r_{L^{\circ }}^{\ast }:\mathcal{H}%
_{s\left( L^{\circ }\right) }\left( V_{2}\right) \overset{\simeq }{%
\rightarrow }\mathcal{H}_{L^{\circ }}\left( V_{1}\right) $.

\begin{lemma}
\label{functor_lemma}For every $\left( M^{\circ },L^{\circ }\right) \in
OLag\left( V_{1}\right) ^{2}$%
\begin{equation*}
T_{M^{\circ },L^{\circ }}\circ r_{L^{\circ }}^{\ast }=r_{M^{\circ }}^{\ast
}\circ T_{s\left( M^{\circ }\right) ,s\left( L^{\circ }\right) }\text{. }
\end{equation*}
\end{lemma}

The system of isomorphisms $\left \{ r_{L^{\circ }}^{\ast }\right \} $
yields the required isomorphism of vector spaces $\mathcal{H}\left( f\right)
:\mathcal{H}\left( V_{2}\right) \overset{\simeq }{\rightarrow }\mathcal{H}%
\left( V_{1}\right) $; evidently, for every $f\in \mathrm{Mor}_{\mathsf{Symp}%
}\left( V_{1},V_{2}\right) $, $g\in \mathrm{Mor}_{\mathsf{Symp}}\left(
V_{2},V_{3}\right) $ 
\begin{equation*}
\mathcal{H}\left( g\circ f\right) =\mathcal{H}\left( g\right) \circ \mathcal{%
H}\left( f\right)
\end{equation*}%
. This concludes the proof of the proposition.

\subsubsection{Proof of Lemma \protect \ref{functor_lemma}}

It will be convenient to work with the system of kernels $\left \{
K_{M^{\circ },L^{\circ }}\right \} $ (Theorem \ref{SS-vN2_prop}). Recall, $%
T_{M^{\circ },L^{\circ }}=I\left[ K_{M^{\circ },L^{\circ }}\right] $. It is
enough to show 
\begin{equation}
r^{\ast }K_{s\left( M^{\circ }\right) ,s\left( L^{\circ }\right)
}=K_{M^{\circ },L^{\circ }},  \label{eq4}
\end{equation}%
for every $\left( M^{\circ },L^{\circ }\right) \in OLag\left( V_{1}\right)
^{2}$. Direct verification reveals that (\ref{eq4}) holds when $\left(
M^{\circ },L^{\circ }\right) \in U_{2}$. Given $\left( M^{\circ },L^{\circ
}\right) \in OLag\left( V_{1}\right) ^{2}$, let $S^{\circ }\in OLag\left(
V_{1}\right) $ be an oriented Lagrangian such that $\left( M^{\circ
},S^{\circ }\right) ,\left( S^{\circ },L^{\circ }\right) \in U_{2}$. Write 
\begin{eqnarray*}
r^{\ast }K_{s\left( M^{\circ }\right) ,s\left( L^{\circ }\right) }
&=&r^{\ast }\left( K_{s\left( M^{\circ }\right) ,s\left( S^{\circ }\right)
}\ast K_{s\left( S^{\circ }\right) ,s\left( L^{\circ }\right) }\right) \\
&=&r^{\ast }K_{s\left( M^{\circ }\right) ,s\left( S^{\circ }\right) }\ast
r^{\ast }K_{s\left( S^{\circ }\right) ,s\left( L^{\circ }\right) } \\
&=&K_{M^{\circ },S^{\circ }}\ast K_{S^{\circ },L^{\circ }}=K_{M^{\circ
},L^{\circ }},
\end{eqnarray*}%
where the first and forth equalities are by the multiplicativity property of
the system $\left \{ K_{M^{\circ },L^{\circ }}\right \} $ (see Theorem \ref%
{SS-vN2_prop}) and the second equality follows from the fact that $r$ is a
morphism of groups.

This concludes the proof of the lemma.

\subsection{Proof of Proposition \protect \ref{Cartesian_prop}}

Let $V_{i}\in \mathsf{Symp}$, $i=1,2$. There are two maps%
\begin{eqnarray*}
r &:&H\left( V_{1}\right) \times H\left( V_{2}\right) \longrightarrow
H\left( V_{1}\times V_{2}\right) , \\
s &:&OLag\left( V_{1}\right) \times OLag\left( V_{1}\right) \longrightarrow
OLag\left( V_{1}\times V_{2}\right) .
\end{eqnarray*}

The map $r$ is a surjective morphism of groups given by $r\left( \left(
v_{1},z_{1}\right) ,\left( v_{2},z_{2}\right) \right) =\left( \left(
v_{1},v_{2}\right) ,z_{1}+z_{2}\right) $. The map $s$ is a bijection of sets
given by 
\begin{equation*}
s\left( L_{1}^{\circ },L_{2}^{\circ }\right) =\left( L_{1}\times
L_{2},o_{L_{1}}\otimes o_{L_{2}}\right) ,
\end{equation*}%
for $L_{i}=\left( L_{i},o_{L_{i}}\right) $, $i=1,2$, where we use the
canonical isomorphism of vector spaces 
\begin{equation*}
\tbigwedge \nolimits^{top}\left( L_{1}\times L_{2}\right) \simeq \tbigwedge
\nolimits^{top}L_{1}\otimes \tbigwedge \nolimits^{top}L_{2}.
\end{equation*}

For every $\left( L_{1}^{\circ },L_{2}^{\circ }\right) \in OLag\left(
V_{1}\right) \times OLag\left( V_{1}\right) $, the map $r$ induces a pull
back isomorphism $r_{\left( L_{1}^{\circ },L_{2}^{\circ }\right) }^{\ast }:%
\mathcal{H}_{s\left( L_{1}^{\circ },L_{2}^{\circ }\right) }\overset{\simeq }{%
\rightarrow }\mathcal{H}_{L_{1}^{\circ }}\otimes \mathcal{H}_{L_{2}^{\circ
}} $.

The proof of the following lemma follows along similar lines as the proof of
Lemma \ref{functor_lemma}.

\begin{lemma}
\label{Cartesian_lemma}For every $\left( L_{1}^{\circ },L_{2}^{\circ
}\right) ,\left( M_{1}^{\circ },M_{2}^{\circ }\right) \in OLag\left(
V_{1}\right) \times OLag\left( V_{1}\right) $%
\begin{equation*}
T_{M_{1}^{\circ },L_{1}^{\circ }}\otimes T_{M_{2}^{\circ },L_{2}^{\circ
}}\circ r_{\left( L_{1}^{\circ },L_{2}^{\circ }\right) }^{\ast }=r_{\left(
M_{1}^{\circ },M_{2}^{\circ }\right) }^{\ast }\circ T_{s\left( M_{1}^{\circ
},M_{2}^{\circ }\right) ,s\left( L_{1}^{\circ },L_{2}^{\circ }\right) }\text{%
. }
\end{equation*}
\end{lemma}

The system $\{r_{\left( L_{1}^{\circ },L_{2}^{\circ }\right) }^{\ast }\}$
gives the required isomorphism $\alpha _{V_{1}\times V_{2}}$.

This concludes the proof of the proposition.

\subsection{Proof of Proposition \protect \ref{duality_prop}}

Let $r:H\left( V\right) \rightarrow H\left( \overline{V}\right) $ be the
isomorphism of Heisenberg groups, given by $r\left( v,z\right) =\left(
v,-z\right) $. For every oriented Lagrangian $L^{\circ }\in OLag\left(
V\right) =OLag\left( \overline{V}\right) $, the map $r$ induces a the
pull-back isomorphism%
\begin{equation*}
r_{L^{\circ }}^{\ast }:\mathcal{H}_{L^{\circ }}\left( \overline{V},\psi
\right) \overset{\simeq }{\rightarrow }\mathcal{H}_{L^{\circ }}\left( V,\psi
^{-1}\right) \text{.}
\end{equation*}

We denote by $\overline{T}_{M^{\circ },L^{\circ }}$ and $T_{M^{\circ
},L^{\circ }}$ the trivializations of $\mathfrak{H}\left( \overline{V},\psi
^{-1}\right) $ and $\mathfrak{H}\left( V,\psi \right) $ respectively.

\begin{lemma}
\label{duality_lemma}For every $\left( M^{\circ },L^{\circ }\right) \in
OLag\left( V\right) ^{2}$%
\begin{equation*}
r_{M^{\circ }}^{\ast }\circ \overline{T}_{M^{\circ },L^{\circ }}=T_{M^{\circ
},L^{\circ }}\circ r_{L^{\circ }}^{\ast }.
\end{equation*}
\end{lemma}

The system $\left \{ r_{L^{\circ }}^{\ast }\right \} $ gives an isomorphism $%
r^{\ast }:\mathcal{H}\left( \overline{V},\psi \right) \overset{\simeq }{%
\rightarrow }\mathcal{H}\left( V,\psi ^{-1}\right) $. Finally, there is a
natural non-degenerate pairing 
\begin{equation*}
\left \langle \cdot ,\cdot \right \rangle _{V}^{\prime }:\mathcal{H}\left(
V,\psi ^{-1}\right) \times \mathcal{H}\left( V,\psi \right) \rightarrow 
%TCIMACRO{\U{2102} }%
%BeginExpansion
\mathbb{C}
%EndExpansion
,
\end{equation*}
induced from the system pairings $\left \langle \cdot ,\cdot \right \rangle
_{L^{\circ }}^{\prime }:\mathcal{H}_{L^{\circ }}\left( V,\psi ^{-1}\right)
\times \mathcal{H}_{L^{\circ }}\left( V,\psi \right) \rightarrow 
%TCIMACRO{\U{2102} }%
%BeginExpansion
\mathbb{C}
%EndExpansion
$, given by 
\begin{equation*}
\left \langle f,g\right \rangle _{L^{\circ }}^{\prime }=\tsum \limits_{v\in
V/L}g\left( v\right) f\left( v\right) ,
\end{equation*}%
for every $f\in \mathcal{H}_{L^{\circ }}\left( V,\psi ^{-1}\right) $ and $%
g\in \mathcal{H}_{L^{\circ }}\left( V,\psi \right) $; we note that the
function $g\cdot f:H\left( V\right) \rightarrow 
%TCIMACRO{\U{2102} }%
%BeginExpansion
\mathbb{C}
%EndExpansion
$ descends to a function on $V/L$. Define 
\begin{equation*}
\left \langle \cdot ,\cdot \right \rangle _{V}=\left \langle r^{\ast }\left(
\cdot \right) ,\cdot \right \rangle _{V}^{\prime }.
\end{equation*}

This concludes the proof of the proposition.

\subsubsection{Proof of Lemma \protect \ref{duality_lemma}}

It is enough to prove the condition for $\left( M^{\circ },L^{\circ }\right)
\in U_{2}$. Let $f\in \mathcal{H}_{L^{\circ }}\left( \overline{V},\psi
\right) $. On the one hand%
\begin{eqnarray*}
r_{M^{\circ }}^{\ast }\overline{T}_{M^{\circ },L^{\circ }}\left[ f\right]
\left( v,z\right) &=&\overline{T}_{M^{\circ },L^{\circ }}\left[ f\right]
\left( v,-z\right) \\
&=&\overline{A}_{M^{\circ },L^{\circ }}\tsum \limits_{m\in M}f\left( m+v,-z-%
\tfrac{1}{2}\omega \left( m,v\right) \right) .
\end{eqnarray*}

On the other hand

\begin{eqnarray*}
T_{M^{\circ },L^{\circ }}\left[ r_{L^{\circ }}^{\ast }\left( f\right) \right]
\left( v,z\right) &=&A_{M^{\circ },L^{\circ }}\tsum \limits_{m\in
M}r_{L^{\circ }}^{\ast }\left( f\right) \left( m+v,z+\tfrac{1}{2}\omega
\left( m,v\right) \right) \\
&=&A_{M^{\circ },L^{\circ }}\tsum \limits_{m\in M}f\left( m+v,-z-\tfrac{1}{2}%
\omega \left( m,v\right) \right) .
\end{eqnarray*}

Noting that $\overline{A}_{M^{\circ },L^{\circ }}=A_{M^{\circ },L^{\circ }}$%
, the result follows.

This concludes the proof of the lemma.

\subsection{Proof of Proposition \protect \ref{reduction_prop}}

Let $V\in \mathsf{Symp}$. Let $I^{\circ }=\left( I,o_{I}\right) $ be an
oriented isotropic subspace in $V$. There is a map $s:OLag\left( I^{\perp
}/I\right) \longrightarrow OLag\left( V\right) $, given by 
\begin{equation*}
s\left( L,o_{L}\right) =\left( pr^{-1}\left( L\right) ,o_{I}\otimes
o_{L}\right) ,
\end{equation*}%
for $\left( L,o_{L}\right) \in OLag\left( I^{\perp }/I\right) $, where $%
pr:I^{\perp }\rightarrow I^{\perp }/I$ is the natural projection and we use
the canonical isomorphism 
\begin{equation*}
\tbigwedge \nolimits^{top}pr^{-1}\left( L\right) \simeq \tbigwedge
\nolimits^{top}I\otimes \tbigwedge \nolimits^{top}L.
\end{equation*}

For every $L^{\circ }\in OLag\left( I^{\perp }/I\right) $ \ there is a
"pull-back" map $r_{L^{\circ }}^{\ast }:\mathcal{H}_{s\left( L^{\circ
}\right) }\left( V\right) \rightarrow \mathcal{H}_{L^{\circ }}\left(
I^{\perp }/I\right) $ defined as follows: Given $f\in \mathcal{H}_{s\left(
L^{\circ }\right) }\left( V\right) $, the function $r_{L^{\circ }}^{\ast
}\left( f\right) $ is the restriction $f_{|Z\cdot I^{\perp }}$, which
descends to a function on $H\left( I^{\perp }/I\right) .$ It is easy to
verify that $r_{L^{\circ }}^{\ast }$ maps $\mathcal{H}_{s\left( L^{\circ
}\right) }\left( V\right) ^{I}$ isomorphically onto $\mathcal{H}_{L^{\circ
}}\left( I^{\perp }/I\right) $.

\begin{lemma}
\label{reduction_lemma}For every $\left( M^{\circ },L^{\circ }\right) \in
OLag\left( I^{\perp }/I\right) $%
\begin{equation*}
F_{M^{\circ },L^{\circ }}\circ r_{L^{\circ }}^{\ast }=r_{M^{\circ }}^{\ast
}\circ F_{s\left( M^{\circ }\right) ,s\left( L^{\circ }\right) }\text{. }
\end{equation*}
\end{lemma}

The system $\{r_{L^{\circ }}^{\ast }\}$ gives the required isomorphism $%
\alpha _{\left( I^{\circ },V\right) }:\mathcal{H}\left( V\right) ^{I}\overset%
{\simeq }{\rightarrow }\mathcal{H}\left( I^{\perp }/I\right) $.

This concludes the proof of the proposition.

\subsubsection{Proof of Lemma \protect \ref{reduction_lemma}}

It will be convenient to work with the canonical system of kernels $\left \{
K_{M^{\circ },L^{\circ }}\right \} $ (Proposition \ref{SS-vN2_prop}), recall
we have $F_{M^{\circ },L^{\circ }}=I\left[ K_{M^{\circ },L^{\circ }}\right] $%
. Let $pr:I^{\perp }\rightarrow I^{\perp }/I$ be the canonical projection
and let $i:I^{\perp }\rightarrow V$ denote the inclusion. It is enough to
show that 
\begin{equation*}
i^{\ast }K_{s\left( M^{\circ }\right) ,s\left( L^{\circ }\right) }=pr^{\ast
}K_{M^{\circ },L^{\circ }}\text{,}
\end{equation*}%
for every $\left( M^{\circ },L^{\circ }\right) \in OLag\left( I^{\perp
}/I\right) ^{2}$. This is done by direct verification using the formulas
from Proposition \ref{explicit_prop}. This concludes the proof of the lemma.

\subsection{Proof of Lemma \protect \ref{isomorphism1_lemma}}

The proof will proceed in several steps.

\textbf{Step} 1. First, we prove that the sheaf $p_{12}^{\ast }\mathcal{%
K\ast }p_{23}^{\ast }\mathcal{K}$ is irreducible $\left[ n_{3}\right] $%
-perverse on $\mathbf{U}_{3}\times \mathbf{H}\left( \mathbf{V}\right) $ and
for this, it is enough to show that $p_{12}^{\ast }\widetilde{\mathcal{K}}%
\mathcal{\ast }p_{23}^{\ast }\widetilde{\mathcal{K}}$ \ is irreducible $%
\left[ n_{3}\right] $-perverse on $\mathbf{U}_{3}\times \mathbf{H}\left( 
\mathbf{V}\right) $.

The convolution $p_{12}^{\ast }\widetilde{\mathcal{K}}\mathcal{\ast }%
p_{23}^{\ast }\widetilde{\mathcal{K}}$ can be written explicitly as\ 
\begin{equation*}
p_{12}^{\ast }\widetilde{\mathcal{K}}\mathcal{\ast }p_{23}^{\ast }\widetilde{%
\mathcal{K}}(N^{\circ },M^{\circ },L^{\circ },h)\simeq \int \limits_{n\in 
\mathbf{N}}\widetilde{\mathcal{K}}\left( N^{\circ },M^{\circ },h\cdot
n^{-1}\right) \otimes \widetilde{\mathcal{K}}(M^{\circ },L^{\circ },n).
\end{equation*}

It is enough to consider the case when $h=m\in \mathbf{M}$ where we get 
\begin{equation}
p_{12}^{\ast }\widetilde{\mathcal{K}}\mathcal{\ast }p_{23}^{\ast }\widetilde{%
\mathcal{K}}(N^{\circ },M^{\circ },L^{\circ },m)\simeq \int \limits_{n\in 
\mathbf{N}}\mathcal{L}_{\psi }\left( \omega \left( n,m\right) \right)
\otimes \widetilde{\mathcal{K}}(M^{\circ },L^{\circ },n).  \label{Fo}
\end{equation}

We see that $p_{12}^{\ast }\widetilde{\mathcal{K}}\mathcal{\ast }%
p_{23}^{\ast }\widetilde{\mathcal{K}}$ $\left( N^{\circ },M^{\circ
},L^{\circ },-\right) _{|\mathbf{M}}$ is an application of $\ell $-adic
Fourier transform to\ $\widetilde{\mathcal{K}}(M^{\circ },L^{\circ },-)_{|%
\mathbf{N}}$, hence, by the Katz--Laumon theorem \cite{KL} $p_{12}^{\ast }%
\widetilde{\mathcal{K}}\mathcal{\ast }p_{23}^{\ast }\widetilde{\mathcal{K}}$
is irreducible $[n_{3}]$-perverse.

\begin{remark}
A more functorial way to describe the previous development is as follows:
Let $\mathbf{C}_{i}\rightarrow \mathbf{U}_{3}$, $i=1,2,3$, denote the
tautological vector bundle with fibers $\left( \mathbf{C}_{i}\right)
_{|\left( L_{1}^{\circ },L_{2}^{\circ },L_{3}^{\circ }\right) }=\mathbf{L}%
_{i}$. What we showed is that 
\begin{equation*}
\left( p_{12}^{\ast }\widetilde{\mathcal{K}}\mathcal{\ast }p_{23}^{\ast }%
\widetilde{\mathcal{K}}\right) _{|\mathbf{C}_{2}}\simeq Four_{\mathbf{C}_{2},%
\mathbf{C}_{1}}\left( p_{23}^{\ast }\widetilde{\mathcal{K}}_{|\mathbf{C}%
_{1}}\right) ,
\end{equation*}%
where $Four_{\mathbf{C}_{2},\mathbf{C}_{1}}$ is the $\ell $-adic Fourier
transform from $\mathsf{D}^{b}\left( \mathbf{C}_{1}\right) $ to $D^{b}\left( 
\mathbf{C}_{2}\right) $, induced from the pairing $w:\mathbf{C}_{1}\times _{%
\mathbf{U}_{3}}\mathbf{C}_{2}\rightarrow \mathbb{G}_{a}{}_{,\mathbf{U}_{3}}$.
\end{remark}

\textbf{Step} 2. It is enough to show that the sheaves $p_{13}^{\ast }%
\mathcal{K}$ and $p_{12}^{\ast }\mathcal{K\ast }p_{23}^{\ast }\mathcal{K}$
are isomorphic on the zero section $\mathbf{U}_{3}\mathbf{\times }\left \{
0\right \} \subset \mathbf{U}_{3}\mathbf{\times H}\left( \mathbf{V}\right) $%
. We have%
\begin{equation}
p_{12}^{\ast }\mathcal{K\ast }p_{23}^{\ast }\mathcal{K}\simeq p_{12}^{\ast }%
\mathcal{A}\otimes p_{23}^{\ast }\mathcal{A}\otimes p_{12}^{\ast }\widetilde{%
\mathcal{K}}\mathcal{\ast }p_{23}^{\ast }\widetilde{\mathcal{K}}.  \notag
\end{equation}

Direct computation reveals that%
\begin{equation*}
\ p_{12}^{\ast }\widetilde{\mathcal{K}}\mathcal{\ast }p_{23}^{\ast }%
\widetilde{\mathcal{K}}(N^{\circ },M^{\circ },L^{\circ },0)\simeq \int
\limits_{n\in \mathbf{N}}\mathcal{L}_{\psi }(\tfrac{1}{2}\omega \left(
r^{M}\left( n\right) ,n\right) ),
\end{equation*}%
where $r^{M}:\mathbf{N}\rightarrow \mathbf{L}$ is the linear map defined by
the relation $r^{M}\left( n\right) -n\in \mathbf{M}$ for every $n\in \mathbf{%
N}$. The map $r^{M}$ induces a map $r_{\wedge }^{M}:\tbigwedge
\nolimits^{top}\mathbf{N\rightarrow }\tbigwedge \nolimits^{top}\mathbf{L}$.
There exists an isomorphism 
\begin{equation*}
\int \limits_{n\in \mathbf{N}}\mathcal{L}_{\psi }\left( \frac{1}{2}\omega
\left( r^{M}\left( n\right) ,n\right) \right) \simeq \mathcal{G}%
_{1}^{\otimes n}\otimes \mathcal{L}_{\sigma }(\left( -1\right) ^{\left( 
\QATOP{n}{2}\right) }\omega _{\wedge }\left( r_{\wedge }^{M}\left(
o_{N}\right) ,o_{N}\right) ).
\end{equation*}

Combining everything together we obtain 
\begin{eqnarray*}
&&p_{12}^{\ast }\mathcal{K\ast }p_{23}^{\ast }\mathcal{K}(N^{\circ
},M^{\circ },L^{\circ },0) \\
&\simeq &\mathcal{G}_{1}^{\otimes 3n}\otimes \mathcal{L}_{\sigma }(\left(
-1\right) ^{\left( \QATOP{n}{2}\right) }\omega _{\wedge }\left(
o_{M},o_{N}\right) \cdot \omega _{\wedge }\left( o_{L},o_{M}\right) \cdot
\omega _{\wedge }\left( r_{\wedge }^{M}\left( o_{N}\right) ,o_{N}\right)
))[4n]\left( 2n\right) \\
&\simeq &\mathcal{G}_{1}^{\otimes n}\otimes \mathcal{L}_{\sigma }(\left(
-1\right) ^{\left( \QATOP{n}{2}\right) +n}\omega _{\wedge }\left(
o_{M},o_{N}\right) \cdot \omega _{\wedge }\left( o_{L},o_{M}\right) \cdot
\omega _{\wedge }\left( r_{\wedge }^{M}\left( o_{N}\right) ,o_{N}\right) ))%
\left[ 2n\right] \left( n\right) \\
&\simeq &\mathcal{G}_{1}^{\otimes n}\otimes \mathcal{L}_{\sigma }(\left(
-1\right) ^{\left( \QATOP{n}{2}\right) }\omega _{\wedge }\left(
o_{L},o_{N}\right) )\left[ 2n\right] \left( n\right) \\
&\simeq &p_{13}^{\ast }\mathcal{K}(N^{\circ },M^{\circ },L^{\circ },0).
\end{eqnarray*}%
where in the second isomorphism we used $\mathcal{G}_{1}^{\otimes 2n}\simeq 
\mathcal{L}_{\sigma }(\left( -1\right) ^{n})[-2n](-n)$ and the third
isomorphism follows from the identity 
\begin{equation*}
\omega _{\wedge }\left( r_{\wedge }^{M}\left( o_{N}\right) ,o_{N}\right)
)=\left( -1\right) ^{n}\omega _{\wedge }\left( o_{N},o_{L}\right) \cdot
\omega _{\wedge }\left( o_{M},o_{N}\right) \omega _{\wedge }\left(
o_{M},o_{L}\right) ^{-1},
\end{equation*}%
which is proved similarly to Lemma \ref{identity_lemma}.

This concludes the proof of Lemma \ref{isomorphism1_lemma}.

\subsection{Proof of Lemma \protect \ref{isomorphism2_lemma}}

Lemma \ref{isomorphism1_lemma} implies that the sheaves $p_{13}^{\ast }%
\mathcal{K}$ and $p_{12}^{\ast }\mathcal{K\ast }p_{23}^{\ast }\mathcal{K}$
are isomorphic on the open subvariety $\mathbf{U}_{3}\times \mathbf{H}\left( 
\mathbf{V}\right) \mathbf{\subset V}_{3}\times \mathbf{H}\left( \mathbf{V}%
\right) $. The sheaf $p_{13}^{\ast }\mathcal{K}$ is irreducible $\left[ n_{3}%
\right] $-perverse as a pull-back by a smooth, surjective with connected
fibers morphism, of an irreducible $[n_{2}]$-perverse sheaf on $\mathbf{OLag}%
\left( \mathbf{V}\right) ^{2}\times \mathbf{H}\left( \mathbf{V}\right) $.
Hence, it is enough to show that the sheaf $p_{12}^{\ast }\mathcal{K\ast }%
p_{23}^{\ast }\mathcal{K}$ is irreducible $\left[ n_{3}\right] $-perverse.

The last assertion follows from the fact that $p_{12}^{\ast }\mathcal{K\ast }%
p_{23}^{\ast }\mathcal{K}$ is principally an application of a, properly
normalized, Fourier transform (see Formula (\ref{Fo})) on $p_{23}^{\ast }%
\mathcal{K}$, hence by the Katz--Laumon theorem \cite{KL} it is irreducible $%
\left[ n_{3}\right] $-perverse.

Let us summarize. We showed that both sheaves $p_{13}^{\ast }\mathcal{K}$
and $p_{12}^{\ast }\mathcal{K\ast }p_{23}^{\ast }\mathcal{K}$ are
irreducible $[n_{3}]$-perverse and are isomorphic on an open subvariety.
This implies that they must be isomorphic. This concludes the proof of the
lemma.

\end{document}